\documentclass[12pt]{amsart}

\usepackage{epsf}

\def\P{{\Bbb P}}
\def\C{{\Bbb C}}
\def\R{{\Bbb R}}

\def\cpt{{\C\P^2}}

\newcommand{\prodlim}{\prod\limits}

\newtheorem{thm}{Theorem}[section]

\newtheorem{prop}[thm]{Proposition}

\newtheorem{lem}[thm]{Lemma}

\newtheorem{rem}[thm]{Remark}

\newtheorem{cor}[thm]{Corollary}

\newtheorem{obs}[thm]{Observation}

\newtheorem{exa}[thm]{Example}

\newtheorem{notation}[thm]{Notation}

\theoremstyle{definition}

\newtheorem{defn}[thm]{Definition}

\newcommand{\een}{\end{enumerate}}

\newcommand{\blem}{\begin{lem}}

\newcommand{\elem}{\end{lem}}

\newcommand{\bcl}{\begin{clm}}

\newcommand{\ecl}{\end{clm}}

\newcommand{\bthm}{\begin{thm}}

\newcommand{\ethm}{\end{thm}}

\newcommand{\bpr}{\begin{prop}}

\newcommand{\epr}{\end{prop}}

\newcommand{\bco}{\begin{cor}}

\newcommand{\eco}{\end{cor}}

\newcommand{\bcon}{\begin{conj}}

\newcommand{\econ}{\end{conj}}

\newcommand{\bde}{\begin{defn}}

\newcommand{\ede}{\end{defn}}

\newcommand{\bex}{\begin{exa}}

\newcommand{\eexa}{\end{exa}}

\newcommand{\bobs}{\begin{obs}}

\newcommand{\eobs}{\end{obs}}

\newcommand{\bexe}{\begin{exe}}

\newcommand{\eexe}{\end{exe}}

\def\prodl{\prod\limits}

\newcommand{\Z}{{\Bbb Z}}

\begin{document}

\pagenumbering{arabic}

\def\qed{\hfill\rule{2mm}{2mm}}

\title[Fundamental groups of tangent conics and tangent lines]{On the fundamental group of the complement of two
tangent conics and an arbitrary number of tangent lines}

\author[Meirav Amram, David Garber and Mina Teicher]{Meirav Amram$^{1,2}$, David Garber$^{2}$ and Mina Teicher}

\stepcounter{footnote} \footnotetext{Partially supported by the Emmy Noether Institute
 Fellowship (by the Minerva Foundation of Germany).}
\stepcounter{footnote} \footnotetext{Partially supported by a
grant from the Ministry of Science, Culture and Sport, Israel and
the Russian Foundation for Basic research, the Russian
Federation.}

\address{Meirav Amram and David Garber, Department of Applied Mathematics,
Faculty of Sciences, Holon Institute of Technology, 52 Golomb St., PO Box 305, 58102
Holon, Israel}

\email{ameirav@math.huji.ac.il,garber@hit.ac.il}

\address{Mina Teicher, Department of Mathematics, Bar-Ilan University, 52900
Ramat-Gan, Israel}

\email{teicher@macs.biu.ac.il}

\date{\today}

\maketitle

\begin{abstract}

We compute  the fundamental groups of the complements of the
family of real conic-line arrangements with up to two conics which
are tangent to each other at two points, with an arbitrary number
of tangent lines to both conics. All the resulting groups turn out
to be big.

\end{abstract}

\maketitle

\medskip


\bigskip

\section{Introduction}\label{intro}

Line arrangements as simply as they are, carry many open questions
around them (including topological and combinatorial questions),
which are slowly solved. In this paper, we go up to conic-line
arrangements, where the parallel questions about them are even
less understood.

We prove here a general theorem on the fundamental group's
structure of complements of families of conic-line arrangements.
The only known results so far in this direction are \cite{AmTe2}
and \cite{AmTeUl}.

In general, the fundamental group of complements of plane curves
is an important topological invariant with many different
applications. Unfortunately, it is hard to compute.

This invariant was used by Chisini \cite{chisini},
Kulikov \cite{Kul} and Kulikov-Teicher \cite{KuTe} in order to
distinguish between different connected components of the moduli space of
surfaces of general type.

Moreover, Zariski-Lefschetz hyperplane section theorem (see
\cite{milnor}) stated that
$$\pi_1 (\P^N - S) \cong \pi_1 (H - H \cap S),$$
where $S$ is an hypersurface and $H$ is a generic 2-plane. Since
$H \cap S$ is a plane curve, we can compute $\pi_1 (\P^N - S)$ in
an easier way, by computing the fundamental group of the complement of
a plane curve.

A different direction for the need of fundamental groups'
computations is for getting more examples of Zariski pairs
(\cite{Z1},\cite{Z2}). A pair of plane curves is called {\it a
Zariski pair} if the two curves have the same combinatorics (i.e. the same
singular points and the same arrangement of irreducible
components), but their complements are not homeomorphic. Some
examples of Zariski pairs can be found at \cite{AB}, \cite{AB-C},
\cite{Deg}, \cite{Ga1}, \cite{oka}, \cite{shimada}, \cite{Tok}
and \cite{uludag}.

\medskip

Let $T_{n,m}$ be the family of real conic-line arrangements in
$\cpt$ with up to two conics, which are tangent to each other at
two points, and with an arbitrary number of tangent lines to each
one of the conics. The main result of this paper is the
computation of the fundamental group $\pi_1(\cpt-T_{n,m})$:

\begin{thm}\label{presentation_Tnm}
Let $Q_1,Q_2$ be two tangent conics in $\cpt$ and let
$\{L_i\}_{i=1}^n$ and $\{L'_j\}_{j=1}^m$ be $n$ and $m$ lines
which are tangent to $Q_2$ and $Q_1$ respectively, see Figure
\ref{Tnm}. Denote:

$$T_{n,m}=Q_1\cup Q_2\cup (\bigcup_{i=1}^n L_i) \cup (\bigcup_{j=1}^m L'_j).$$

\begin{figure}[h]
\epsfysize=6cm \centerline{\epsfbox{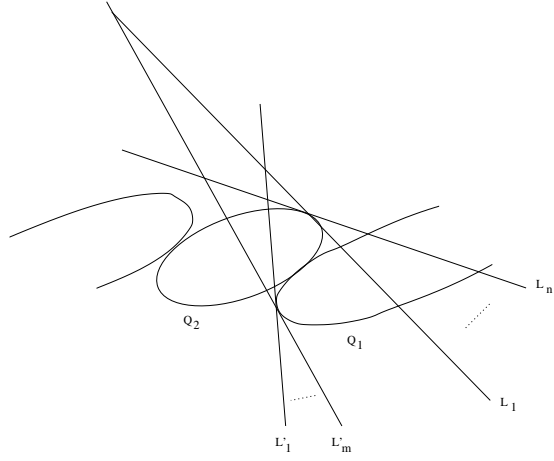}}
\caption{The arrangement $T_{n,m}$}\label{Tnm}
\end{figure}

Then the fundamental group $\pi_1(\cpt-T_{n,m})$ is generated by the generators
$x_{2}$ (related to $Q_1$), $x_{5}$ (related to $Q_2$), $x_{6},\dots, x_{n+4}$ (related to $L_2,\dots, L_{n}$),
 $x_{n+5},\dots, x_{n+m+4}$ (related to $L'_1,\dots, L'_m$). 
This group admits the following relations:

\begin{enumerate}

\item $x_{n+5}^{-1} x_5 x_{n+5} = x_{n+6}^{-1} \cdots  x_{n+m+4}^{-1} \cdot x_{5} \cdot x_{n+m+4} \cdots  x_{n+6}$

\item $(x_{2} x_i)^2 = (x_{i} x_{2})^2$, where $i=5, 6, \dots, n+4$

\item $(x_{5} x_i)^2 = (x_{i} x_5)^2$, where $n+5 \leq i \leq n+m+4$

\item $[x_i, x_j]=e$, where $6 \leq i \leq n+4$ and $j=5, n+5, \dots, n+m+4$

\item $[x_i, x_2 x_{5} x_2^{-1}]$, where $6 \leq i \leq n+4$

\item $[x_{2}^{-1} x_i x_{2}, x_{j}]=e$, where $6 \leq i < j \leq n+4$

\item $[x_2, x_{i}]=e$, where  $n+6 \leq i \leq n+m+4$

\item $[x_{5} x_{2} x_{5}^{-1}, x_i]=e$, where $n+5 \leq i \leq n+m+4$

\item $[x_{n+5}, x_{i}]=e$, where $i= 2, n+6, \dots, n+m+4$
\item $[x_{5}^{-1} x_{i} x_{5}, x_j]=e$, where $n+6 \leq i < j \leq n+m+4$.
\end{enumerate}

\end{thm}

We exclude here the case where at least one line is tangent to both conics 
(in their tangency point).

These arrangements may appear as a branch curve of a generic projection
to $\cpt$ of a surface of general type (see for example \cite{Hi}).

\medskip

A group is called {\it big}  if it contains a subgroup,
generated by two or more generators, which is free. By the above
computations, we have the following corollary:

\begin{cor}\label{bigness}

The fundamental group $\pi_1(\cpt-T_{n,m})$ is big.

\end{cor}

Algorithmically, this paper uses local computations (local braid
monodromies and their induced relations, see also \cite{AmGaTe}),
braid monodromy techniques of Moishezon-Teicher (see \cite{GaTe},
\cite{Mo1}, \cite{MoTe1}, \cite{MoTe2}, \cite{MoTe3} and
\cite{MoTe4}), the van Kampen Theorem (see \cite{vK}) and some
group simplifications and calculations for studying the
fundamental groups (see also \cite{AmOg} and \cite{AmTe}). 

Note that the new tool we are using in this paper is to 
construct the BMF of the general curve in an inductive
way.

\medskip

The paper is organized as follows. Section \ref{sec2} deals with
the fundamental group of a family of curves consisting of a conic
and $n$ tangent lines. In Section \ref{sec3}, we deal with the
fundamental group of a family of curves consisting of two tangent
conics and $n$ tangent lines to one of the conics. In Section
\ref{sec4}, we compute the most general case, where we have two
tangent conics and two sets of lines, each of which is tangent to
one of the conics.

In each section, we start with the computation
of the braid monodromy of the curves by Moishezon-Teicher
techniques and then we compute the simplified presentation of the
fundamental groups by the van Kampen Theorem and group theoretic
calculations.

\section{The fundamental group of the complement of a conic with $n$ tangent
 lines}\label{sec2}

In this section, we compute the fundamental group of the
complement of a conic and $n$ tangent lines. We start by
illustrating the cases of one and two lines which are tangent to a conic,
and then we present the computation of the general case of a conic
and $n$ tangent lines.

Note that this family of conic-line arrangements was already
computed in \cite{AmTeUl} too (the family is denoted there by
$A_i$). Although the presentations given here are different from
those appeared in \cite{AmTeUl}, it can be easily shown that the
two presentations yield isomorphic groups.

We recompute their presentations here (by different methods) for
the sake of completeness of presentations of conic-line
arrangements with up to two tangent conics. One more reason for
this is that it illustrates the method of computing the braid
monodromy factoriazations  in an inductive way, as we will do in the case of two
tangent conics and any number of tangent lines.

\medskip

In the first subsection, we present the braid monodromy
factorizations of this type of curves and in the next subsection,
we compute the corresponding fundamental groups.

\subsection{Braid monodromy factorizations}\label{bmfcn}

Here we compute the braid monodromy factorizations (BMFs, see
\cite{KuTe}) which correspond to a conic with an arbitrary number
of lines tangent to it. We start with some notations.

\begin{notation}\label{nodes}
Let $z_{ij}$ (resp. $\bar{z}_{ij}$) be a path below (resp. above) the axis, which connects points $i$ and
$j$.
We denote by $Z_{i j}$ (resp. $\bar{Z}_{ij}$) the counterclockwise half-twist of $i$ and $j$ along
$z_{ij}$ (resp. $\bar{z}_{ij}$).

We denote $Z^2_{i,j j'} = Z_{i j'}^2 \cdot Z_{i  j}^2$ and $Z^2_{i i', j  j'} = Z_{i', j  j'}^2
\cdot Z_{i, j  j'}^2$.

Conjugation of braids is denoted by $a^b = b^{-1}ab$.
\end{notation}

\begin{exa}\label{example6}
The path in Figure \ref{example}(a) is constructed as follows:
take a path $z_{34}$ and conjugate it by the full-twist
$\bar{Z}_{13}^2$. Its corresponding half-twist is
$Z_{34}^{\bar{Z}_{13}^2}$.
 The path in Figure \ref{example}(b) is constructed as follows:
 take again $z_{34}$ and conjugate it first by $Z_{23}^2$  and then by $Z_{13}^2$.
 It corresponding half-twist is $Z_{34}^{Z_{23}^2 Z_{13}^2}$.
\begin{figure}[h!]
\epsfxsize=7cm 
\begin{minipage}{\textwidth}
\begin{center}
\epsfbox{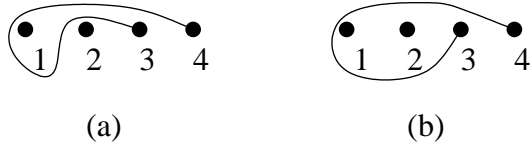}
\end{center}
\end{minipage}
\caption{Examples of conjugated braids}\label{example}
\end{figure}
\end{exa}

Here, we start the computation of the BMFs.

\begin{lem}\label{1conic_2+1}
Let $Q$ be a conic in $\cpt$ and let $L$ be a line which is
tangent to $Q$ (see Figure \ref{gen_cl1}). Let $C_1 =Q \cup L$.

\begin{figure}[h]
\epsfysize=3cm
\centerline{\epsfbox{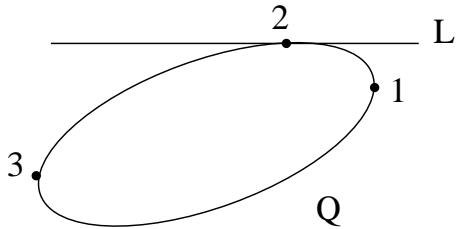}}
\caption{The arrangement $C_1$}\label{gen_cl1}
\end{figure}

Then its BMF is  $\Delta_{C_1}^2 = F_1 \cdot (Z_{1 1'})^{Z_{1' 2}^{-2}}$,
where $F_1 = Z_{1 1'} \cdot (Z_{12}^4)^{Z_{1 1'}^2}$.
\end{lem}

\begin{proof}
The computation of the BMF is done in a local model as follows: we
start the computation by taking a typical fiber and enumerate the
points of the curve in this fiber from $1$ to $n$ ($n$ is the
degree of the given curve). At this stage, we classify the
singularities according to their types, and we summarize it in a table:

$$
\begin{array}{|c||c|c|c|}
\hline
1 & P_1                 & 1 & \Delta ^{1\over 2}_{\R I_2} \langle 1 \rangle \\
2 & \langle 2,3 \rangle & 4 & \Delta ^2 \langle 2,3 \rangle \\
3 & \langle 1,2 \rangle & 1 & \Delta ^{1\over 2}_{I_2 \R} \langle 1 \rangle \\
\hline
\end{array}
$$

The second column contains the Lefschetz pairs corresponding to
the singularities. In the third column, we indicate the types of
the singularities, and in the last column we have the
corresponding local diffeomorphisms of the singularities (see
\cite{AmTe}).

In the next step, we apply the braid monodromy technique of
Moishe\-zon-Teicher \cite{MoTe2} (see also \cite{AmGaTe} and
\cite{AmTe}). In Figure \ref{gen_cl1_bm}, we show the skeletons of
the braids, obtained by the technique. Note that we skip some
steps in the braid monodromy's computation. Note also that at the
end of the computation, we enumerate again the points in the fiber
according to the global settings (for example, the two points of a
conic are denoted by $i$ and $i'$).

\begin{figure}[ht!]
\epsfysize=4.5cm \centerline{\epsfbox{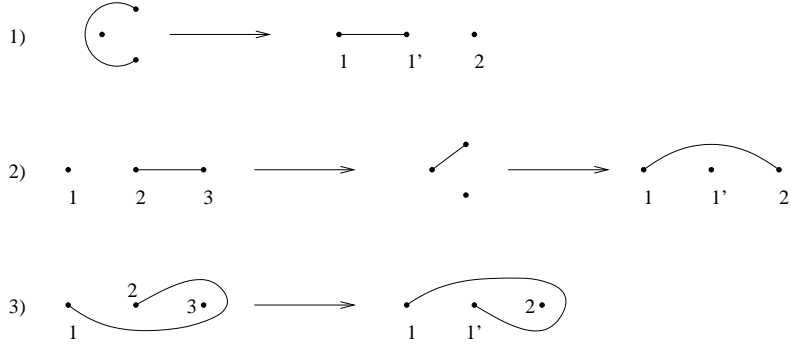}}
\caption{The braid monodromy related to $C_1$}\label{gen_cl1_bm}
\end{figure}

The first braid is $Z_{1 1'}$, the second one is  $(Z_{12}^4)^{Z_{1 1'}^2}$ and
the third one is $(Z_{1 1'})^{Z_{1' 2}^{-2}}$, and hence the braid
monodromy factorization is as stated.
\end{proof}

\begin{lem}\label{1conic_2+2}
Let $Q$ be a conic in $\cpt$ and let $L_1,L_2$ be two lines which are
tangent to $Q$ (see Figure \ref{gen_cl2}). Let $C_2 = Q \cup L_1
\cup L_2$.

\begin{figure}[ht!]
\epsfysize=3cm
\centerline{\epsfbox{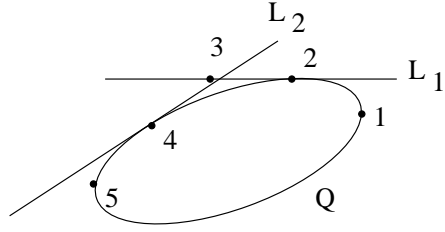}}
\caption{The arrangement $C_2$}\label{gen_cl2}
\end{figure}

Then its BMF is $\Delta_{C_2}^2 = F_1 \cdot F_2 \cdot (Z_{1 1'})^{Z_{1' 2}^{-2} Z_{1' 3}^{-2}}$,
where $F_1 = Z_{1 1'} \cdot (Z_{12}^4)^{Z^2_{1 1'}}$ and $F_2 = (Z_{23}^2)^{\bar{Z}_{1 2}^{2}} \cdot
\bar{Z}_{13}^4$.

\end{lem}

\begin{proof}
The table for computing the braid monodromy is:
$$
\begin{array}{|c||c|c|c|}
\hline
1 & P_1                 & 1 & \Delta ^{1\over 2}_{\R I_2} \langle 1 \rangle \\
2 & \langle 2,3 \rangle & 4 & \Delta ^2 \langle 2,3 \rangle \\
3 & \langle 3,4 \rangle & 2 & \Delta \langle 3,4 \rangle \\
4 & \langle 2,3 \rangle & 4 & \Delta ^2 \langle 2,3 \rangle \\
5 & \langle 1,2 \rangle & 1 & \Delta ^{1\over 2}_{I_2 \R} \langle 1 \rangle
\\
\hline
\end{array}
$$

In a similar way as in the previous lemma, we get the skeletons for the braid
monodromies presented in Figure \ref{gen_cl2_bm}.

\begin{figure}[ht]
\epsfysize=7cm \centerline{\epsfbox{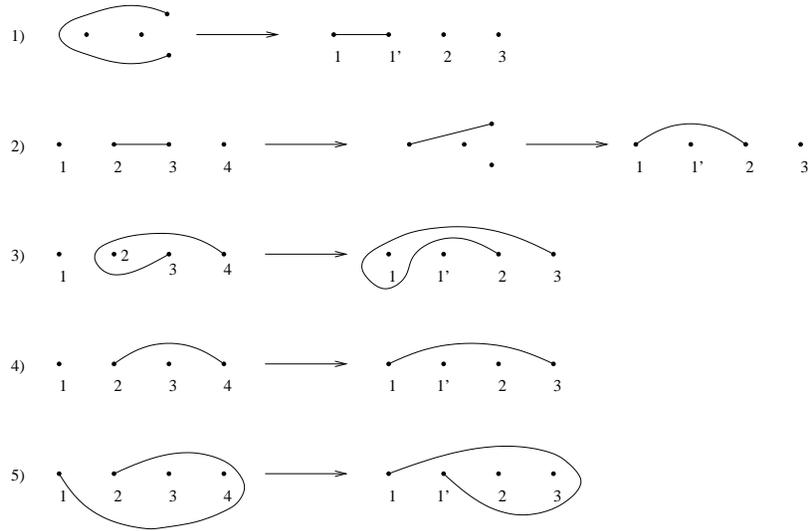}}
\caption{The braid monodromy related to $C_2$ }\label{gen_cl2_bm}
\end{figure}
\end{proof}

\begin{thm}\label{thmbmf-cn}
Let $Q$ be a conic in $\cpt$ and let $L_1,\dots, L_n$ be $n$ lines which are
tangent to $Q$. Let $C_n=Q\cup \bigcup_{i=1}^n L_i$ (see Figure
\ref{C_n}).

\begin{figure}[h]
\epsfysize=6cm \centerline{\epsfbox{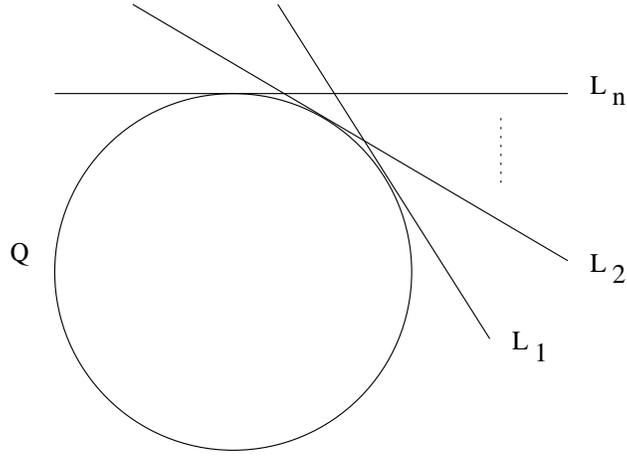}}
\caption{The arrangement $C_n$}\label{C_n}
\end{figure}

Then its BMF is:
$$\Delta_{C_n}^2 = F_1 \cdot F_2 \cdots F_n \cdot (Z_{1 1'})^{\displaystyle \prod_{i=n+1}^2 Z_{1',i}^{-2}},$$
where:
$$
F_1 = Z_{1 1'} \cdot (Z_{12}^4)^{Z_{1 1'}^2}
$$
and
$$F_i = \prod_{k=2}^i (Z_{k,i+1}^2)^{\bar{Z}_{1k}^{-2}} \cdot \bar{Z}_{1,i+1}^4, \quad 2 \leq i \leq n.$$
\end{thm}

\begin{proof}
We prove it by induction. The cases $n=1,2$ have been proved in
the previous lemmas.

Assume that $\Delta_{C_n}^2$ is as given above. We will show that:
$$\Delta_{C_{n+1}}^2 = F_1 \cdot F_2 \cdots F_n \cdot F_{n+1} \cdot (Z_{1
1'})^{\displaystyle \prod_{i=n+2}^2 Z_{1',i}^{-2}}.$$

The sub-factorization $F_1 \cdot F_2 \cdots F_n$ appears in
$\Delta_{C_{n+1}}^2$, since by the structure of the configuration, 
 the additional line $L_{n+1}$
is located to the left of the rest of the lines $L_1, \dots, L_n$.

Note that by induction, adding the line $L_{n+1}$, causes no change in the
factors $F_i$s, $1 \leq i \leq n$. This is due to the appearance
of this line on the left side of $L_1, \dots, L_n$.

Now we explain the term:
$$F_{n+1} = {\Big(}\prod_{k=2}^{n+1} (Z_{k,n+2}^2)^{\bar{Z}_{1k}^{-2}}{\Big)} \cdot \bar{Z}_{1,n+2}^4.$$
The first part  $\prod (Z_{k,n+2}^2)^{\bar{Z}_{1k}^{-2}}$ is the
product of braids which correspond to the intersection points of
$L_{n+1}$ with the other lines $L_1,\dots, L_n$. Now,
$\bar{Z}_{1,n+2}^4$ is the braid which corresponds to the tangency
point of the line $L_{n+1}$ and the conic $Q$.

We still have to explain the conjugated element. We have:
$$(Z_{11'})^{\tiny \prod_{i=n+2}^2 Z_{1',i}^{-2}} = \left( (Z_{1 1'})^{\tiny
\prod_{i=n+1}^2 Z_{1',i}^{-2}} \right) ^{Z_{1',n+2}^{-2}},$$ where
$(Z_{1 1'})^{\tiny \prod_{i=n+1}^2 Z_{1',i}^{-2}}$ appears in
$\Delta^2_{C_n}$. Due to the addition of the line $L_{n+1}$, we have to
conjugate this braid by the braid  $Z_{1',n+2}^{-2}$.
\end{proof}

\subsection{The corresponding fundamental groups}\label{fgcn}
Based on the BMFs we computed in the previous subsection, we now
compute the fundamental groups of the affine and projective
complements of the curve $C_n$.

\begin{lem}
Let $C_1$ be as in Lemma \ref{1conic_2+1} (see Figure
\ref{gen_cl1}).

Then we have:
$$\pi_1(\C^2-C_1) \cong \Z ^2 \qquad ; \qquad \pi_1(\cpt-C_1) \cong \Z.$$
\end{lem}

\begin{proof}
By applying the van Kampen Theorem \cite{vK} on $\Delta_{C_1}^2$
(computed in  Lemma \ref{1conic_2+1}), we get the following presentation
for $\pi_1(\C^2-C_1)$:\\
Generators: $\{ x_1,x_{1'},x_2 \}$ \\
Relations:
\begin{enumerate}
\item $x_1= x_{1'}$
\item $(x_{1'} x_1 x_{1'}^{-1} x_2)^2= (x_2 x_{1'} x_1 x_{1'}^{-1})^2$
\item $x_1 = x_{1'}^{-1} x_2^{-1} x_{1'} x_{2} x_{1'}$
\end{enumerate}

By Relation (1), the generator $x_{1'}$ is redundant, and by
Relations (1) and (3), we get $x_1 x_2 = x_2 x_1$. Hence, we have:
$$\pi_1(\C^2-C_1) = \langle x_1,x_2\ |\ x_1 x_2 = x_2 x_1
\rangle \cong \Z^2.$$

Adding the projective relation $x_2 x_{1'} x_1 =e$, we easily get:
$$\pi_1(\cpt-C_1) \cong \Z.$$
\end{proof}

\begin{lem}
Let $C_2$ be as in Lemma \ref{1conic_2+2} (see Figure
\ref{gen_cl2}).

Then:
$$\pi_1(\C^2-C_2) \cong \left\langle x_1,x_2,x_3 \left| \begin{array}{c}
(x_1 x_2)^2= (x_2 x_1)^2 \\ (x_1 x_3)^2= (x_3 x_1)^2 \\
\ [x_3,x_1^{-1} x_2 x_1]=[x_3 x_2,x_1]=e \end{array} \right.
\right\rangle, $$ where $[a,b]=aba^{-1}b^{-1}$ is the commutator
element of $a$ and $b$, and
$$\qquad \pi_1(\cpt-C_2) \cong \langle x_1, x_2 \ |\ (x_1 x_2)^2= (x_2 x_1)^2 \rangle,$$
where  the generator $x_1$ corresponds to the conic and the
generators $x_2,x_3$ correspond to the lines.
\end{lem}

\begin{proof}
By applying the van Kampen Theorem on $\Delta_{C_2}^2$ (computed in Lemma
\ref{1conic_2+2}), we get the following presentation
for $\pi_1(\C^2-C_2)$:\\
Generators: $\{ x_1,x_{1'},x_2,x_3 \}$ \\
Relations:
\begin{enumerate}
\item $x_1= x_{1'}$
\item $(x_{1'} x_1 x_{1'}^{-1} x_2)^2= (x_2 x_{1'} x_1 x_{1'}^{-1})^2$
\item $[x_3, x_2 x_{1'} x_{1} x_{1'}^{-1} x_2 x_{1'} x_{1}^{-1} x_{1'}^{-1} x_2^{-1}]=e$
\item $(x_2 x_{1'} x_{1} x_{1'}^{-1} x_2^{-1} x_3)^2= (x_3 x_2 x_{1'} x_{1} x_{1'}^{-1} x_2^{-1})^2$
\item $x_1 = x_{1'}^{-1} x_2^{-1} x_{3}^{-1} x_{1'} x_3 x_2 x_{1'}$
\end{enumerate}

By Relation (1), the generator $x_{1'}$ is redundant. So we get:\\
Generators: $\{ x_1,x_2,x_3 \}$ \\
Relations:
\begin{enumerate}
\item $(x_1 x_2)^2= (x_2 x_{1})^2$
\item $[x_3, x_2 x_{1} x_2 x_{1}^{-1} x_2^{-1}]=e$
\item $(x_2  x_{1}  x_2^{-1} x_3)^2= (x_3 x_2  x_{1}  x_2^{-1})^2$
\item $x_1 =  x_2^{-1} x_{3}^{-1} x_{1} x_3 x_2$
\end{enumerate}

By Relation (1), Relation (2) is simplified to $[x_3, x_1^{-1} x_2 x_{1}]=e$. Relation (4) is rewritten as
$[x_1, x_3 x_2]=e$. We treat now Relation (3):
$$x_2  x_{1}  x_2^{-1} x_3 x_2  x_{1}  x_2^{-1} x_3 = x_3 x_2  x_{1}  x_2^{-1} x_3 x_2  x_{1}  x_2^{-1}.$$
By the relation  $[x_3, x_1^{-1} x_2 x_{1}]=e$, it can be written as:
$$x_2  x_{1}  x_2^{-1} x_1  x_{3} x_3 = x_1 x_{3}  x_1 x_3.$$
By Relation (4), $x_2  x_{1}  x_2^{-1} = x_{3}^{-1}  x_1 x_3$.
So,  we get $(x_1 x_3)^2=(x_3 x_1)^2$.

The addition of the projective relation $x_3 x_2 x_{1'} x_1 =e$
enables us to omit the generator $x_3=x_1^{-2}x_2^{-1}$. By
substituting $x_3$ in any place it appears, we get the needed
presentation for $\pi_1(\cpt-C_2)$.
\end{proof}

Now, we compute the presentation for the general case.

\begin{prop}
Let $Q$ be a conic in $\cpt$ and let $L_i$ ($1 \leq i \leq n$) be
$n$ lines which are tangent to $Q$ (see Figure \ref{C_n}). Let $C_n = Q
\cup \bigcup_{i=1}^n L_i$.

Then:
$$\pi_1(\C^2-C_n) \cong \left\langle x_1,x_2,\dots, x_{n+1} \left| \begin{array}{c}
(x_1 x_i)^2= (x_i x_1)^2,\quad 2 \leq i \leq n+1 \\
\ [x_j, x_1^{-1} x_i x_1]=e, \quad 2 \leq i < j \leq n+1  \\
\ [x_{n+1} \cdots x_3 x_2,x_1]=e
\end{array} \right. \right\rangle$$
and
$$\pi_1(\cpt-C_n) \cong \left\langle x_1,x_2,\dots, x_{n+1} \left| \begin{array}{c}
(x_1 x_i)^2= (x_i x_1)^2,\quad 2 \leq i \leq n+1 \\
\ [x_j, x_1^{-1} x_i x_1]=e, \quad 2 \leq i < j \leq n+1  \\
x_{n+1} \cdots x_3 x_2 x_1^2=e
\end{array} \right. \right\rangle,$$
where $x_1$ is a generator which corresponds to the conic and $x_2, \dots,
x_{n+1}$ are  generators which correspond to the $n$ lines.
\end{prop}

\begin{proof}
In this case, the set of generators for the presentation of the
fundamental group is  $\{ x_1,x_{1'},x_2,\dots, x_{n+1} \}$. By
applying the van Kampen Theorem on the braid monodromy computed in Theorem \ref{thmbmf-cn},
we get three types of relations:
\begin{enumerate}
\item From the branch points, we get $x_1= x_{1'}$ and
$$x_1 = x_{1'}^{-1} x_2^{-1} x_3^{-1} \cdots x_{n+1}^{-1} x_{1'} x_{n+1}  \cdots x_3 x_2 x_1.$$

\item From the nodal points, we get the following set of relations:
$$[x_j,x_{j-1} x_{j-2} \cdots x_{1'} x_1 x_{1'}^{-1} x_{2}^{-1} \cdots x_{i-1}^{-1} x_{i} x_{i-1}
\cdots x_{2} x_{1'} x_{1}^{-1} x_{1'}^{-1} \cdots x_{j-2}^{-1}
x_{j-1}^{-1}]=e,$$
where $2 \leq i < j \leq n+1$.

 \item From the tangency points, we get: 
$$(x_1 x_{2}^{-1} x_{3}^{-1} \cdots x_{i-1}^{-1} x_{i} x_{i-1} \cdots x_3 x_2)^2 =
(x_{2}^{-1} x_{3}^{-1} \cdots x_{i-1}^{-1} x_{i} x_{i-1} \cdots
x_3 x_2 x_1)^2,$$
where $2 \leq i \leq n+1$.
\end{enumerate}

By the relations of the first type (from the branch points),
$x_{1'}$ is redundant, and we have: $[x_{n+1} \cdots x_3 x_2,x_1]=e$.

Now, we simplify the relations of types (2) and (3). In type (3),
we get for $i=2$: $(x_1 x_{2})^2 = (x_{2} x_{1})^2$. In type (2),
for $i=2$ and $j=3$, we get:
 $[x_3, x_2 x_{1} x_2 x_{1}^{-1} x_2^{-1}]=e$. By $(x_1 x_{2})^2 = (x_{2} x_{1})^2$, it is
 $[x_3, x_{1}^{-1} x_2 x_{1}]=e$. 

Returning again to type (3), for $i=3$, we have:
 $(x_1 x_{2}^{-1} x_{3} x_{2})^2 =  (x_{2}^{-1} x_{3} x_2 x_1)^2$, which is:
$$
(x_1 x_{2}^{-1} x_{3} x_{2}) (x_1 x_{2}^{-1} x_{3} x_{2}) =
(x_{2}^{-1} x_{3} x_2 x_1) (x_{2}^{-1} x_{3} x_2 x_1).
$$
Since  $(x_1 x_{2})^2 = (x_{2} x_{1})^2$, we have
$$
x_1 x_{2}^{-1} x_3 x_1^{-1} x_{2}^{-1} x_1 x_{2} x_1 x_{3} x_{2} = x_{2}^{-1} x_{3} x_1^{-1} x_2^{-1} x_1 x_{2} x_{1} x_{3} x_2 x_1,
$$
and by $[x_3, x_{1}^{-1} x_2 x_{1}]=e$, we get:
$$
x_1 x_{2}^{-1} x_1^{-1} x_{2}^{-1} x_1 x_3 x_{2} x_1 x_{3} x_{2} = x_{2}^{-1} x_1^{-1} x_2^{-1} x_1 x_3 x_{2} x_{1} x_{3} x_2 x_1.
$$
This relation is simplified to (by $(x_1 x_{2})^2 = (x_{2} x_{1})^2$):
$$
x_1 x_3 x_{2} x_1 x_{3} x_{2} = x_3 x_{2} x_{1} x_{3} x_2 x_1.
$$
This relation is:
$$
x_1 x_3 x_{2} x_1 x_{3} = x_3 x_{2} x_{1} x_{3} x_2 x_1 x_2^{-1}.
$$
Using  $(x_1 x_{2})^2 = (x_{2} x_{1})^2$, we have:
$$
x_1 x_3 x_{2} x_1 x_{3} = x_3 x_{2} x_{1} x_{3} x_1^{-1} x_2^{-1} x_1 x_2 x_1.
$$
Now we add $x_1 x_1^{-1}=e$:
$$
x_1 x_3 (x_1 x_1^{-1}) x_{2} x_1 x_{3} = x_3 x_{2} x_{1} x_{3}
x_1^{-1} x_2^{-1} x_1 x_2 x_1.
$$
Since  $[x_3, x_{1}^{-1} x_2 x_{1}]=e$,
$$
x_1 x_3 x_1 x_3 x_1^{-1} x_{2} x_1 = x_3 x_{2} x_{1} x_1^{-1} x_2^{-1} x_1 x_3 x_2 x_1.
$$
And we get $(x_1 x_{3})^2 = (x_{3} x_{1})^2$.

Now we return again to type (2), for $i=2$ and $j=4$, we have:
$[x_4, x_3 x_{2} x_1 x_{2} x_{1}^{-1} x_2^{-1} x_{3}^{-1}]=e$. By
$(x_1 x_{2})^2 = (x_{2} x_{1})^2$, it is: $[x_4, x_3 x_1^{-1}
x_{2} x_{1} x_{3}^{-1}]=e$, which is simplified to: $[x_4,
x_1^{-1} x_{2} x_{1}]=e$, by using $[x_3, x_1^{-1} x_2 x_1]=e$. 

Now, for $i=3$ and $j=4$, we have: $[x_4, x_3 x_2 x_1 x_2^{-1}
x_{3} x_{2} x_{1}^{-1} x_{2}^{-1} x_{3}^{-1}]=e$. By $(x_1
x_{2})^2 = (x_{2} x_{1})^2$, we rewrite it as:
$$[x_4, x_3 x_1^{-1} x_2^{-1} x_1 x_2 x_{1} x_{3} x_{1}^{-1} x_{2}^{-1} x_{1}^{-1} x_2 x_1 x_3^{-1}]=e,$$
which is (by $[x_3, x_{1}^{-1} x_2 x_{1}]=e$):
$$[x_4, x_1^{-1} x_2^{-1} x_1 x_3 x_2 x_{1} x_{3} x_{1}^{-1} x_{2}^{-1} x_3^{-1} x_{1}^{-1} x_2 x_1]=e.$$
Now, by $[x_4, x_1^{-1} x_{2} x_{1}]=e$:
$$[x_4, x_3 x_2 x_{1} x_{3} x_{1}^{-1} x_{2}^{-1} x_3^{-1}]=e.$$
By  $[x_3, x_{1}^{-1} x_2 x_{1}]=e$, we get:  $[x_4, x_3 x_1 x_3
x_1^{-1} x_3^{-1}]=e$ and by $(x_1 x_{3})^2 = (x_{3} x_{1})^2$, we
get $[x_4, x_1^{-1} x_3 x_1]=e$. 

In the same way, we get the following sets of relations:
\begin{equation}
[x_{j}, x_1^{-1} x_i x_1]=e, \ \ \mbox{where $2 \leq i < j \leq
n+1$}
\end{equation}
and
\begin{equation}
(x_1 x_{i})^2 = (x_{i}x_{1})^2 \ \ \mbox{where \ \ $2 \leq i \leq
n+1$},
\end{equation}
as required.

Adding the projective relation $x_{n+1} x_n \cdots x_2
x_1^2=e$, the relation $$[x_{n+1} \cdots x_3 x_2,x_1]=e$$
is redundant, so we get the needed presentation for the projective
case too.
\end{proof}

\subsection{The groups are big}

Following the results of the previous subsections, we have:
\begin{cor}\label{big}
The affine and projective fundamental groups of a conic-line arrangement,
composed of one conic and at least two tangent lines, are big.
\end{cor}

We have already proved a similar corollary in \cite{AmGaTe}, and
for the sake of completeness, we recall its proof shortly here.

\begin{proof}
Observe that if a group has a big quotient, then the group itself is big.

Denote $G = \pi_1(\cpt-C_2) \cong \langle a, b \ |\ (a b)^2= (b a)^2 \rangle$.
We show that $G$ is big. We start by taking the quotient by the relation $(a
b)^2=e$. Then we get the group $\langle a, b \ |\ (a b)^2= (b a)^2 = e \rangle$.

Let us take new generators $x=ab, y=b$, then the relation $(ab)^2
= (ba)^2 =e$
becomes: $x^2= (yxy^{-1})^2=e$, which is equal to: $x^2= y x^2
y^{-1}=e$. Hence, we have:
$$
\langle a, b \ |\ (a b)^2= (b a)^2 = e \rangle \cong \langle x,y\ |\ x^2=e \rangle \cong \Z * \Z/2,
$$
where $*$ is the free product.

Now, the quotient of $G/$$\langle x^2= e \rangle$ by the subgroup generated by $y^3$ is
$\Z/2 * \Z/3$, which is known to be big. By the observation,
$G = \pi_1(\cpt-C_2)$ is big too.

Since $C_2$ is a sub-arrangement of $C_n$, $n \geq 2$,  the group
$G = \pi_1(\cpt-C_2)$ is a quotient
of the group $\pi_1(\cpt-C_n)$ (by sending  the generators
which correspond to the additional lines to $e$).
Now, since $G$ is big, the above groups are big as well.

The above proof shows that the projective fundamental group of $C_n$, $n \ge 2$, is big.
Since the projective fundamental group is a quotient of the affine fundamental group
(by the projective relation),
the corresponding affine fundamental groups are also big, 
again by the observation, as needed.
\end{proof}

\section{The fundamental group of the complement of two tangent conics with
$n$ tangent lines to the same conic}\label{sec3}

In this section, we compute the fundamental group of the
projective complement of a curve composed of two tangent conics
and $n$ additional lines, which are tangent to the same conic.

As before, in the first subsection, we present the BMF of this
type of curves and in the second subsection, we compute the
corresponding fundamental groups.

\subsection{The braid monodromy factorizations}

Let $T_{n,0}$ be a conic-line arrangement, composed of two tangent
conics (which are tangent to each other at two points), and with $n$ additional lines (which are tangent to the same conic).

We start with the arrangements $T_{0,0}, T_{1,0}$ and $T_{2,0}$.
Their fundamental groups have been already computed in
\cite{AmGaTe}.

We give here their BMFs:
\begin{lem}\label{t00}
Let $T_{0,0}$ be a conic arrangement composed of two  conics which are tangent at two points (see Figure \ref{fig_t00}).
\begin{figure}[!ht]
\epsfysize=3cm
\centerline{\epsfbox{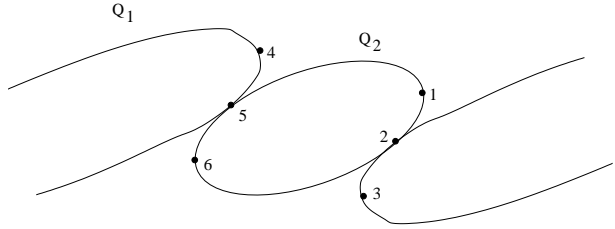}}
\caption{The arrangement $T_{0,0}$}\label{fig_t00}
\end{figure}

The BMF of $T_{0,0}$ is:
$$\Delta^2_{T_{0,0}} = Z_{3 4} \cdot Z^4_{2 4} \cdot (Z_{1 2})^{Z^2_{2 4}}
 \cdot
(Z_{1 2})^{Z^2_{2 3}} \cdot Z^4_{2 3} \cdot  Z_{3 4}.$$
\end{lem}

\begin{lem}\label{t10}
Let $T_{1,0}$ be a conic-line arrangement composed of two  conics which are tangent at
 two points, and a line which is tangent to one of the conics (see Figure \ref{fig_t10}).
\begin{figure}[!ht]
\epsfysize=3.5cm
\centerline{\epsfbox{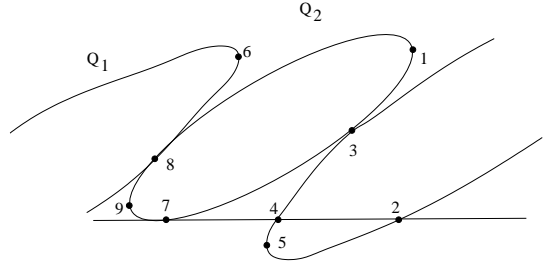}}
\caption{The arrangement $T_{1,0}$}\label{fig_t10}
\end{figure}

Then its BMF is:
\begin{eqnarray*}
\Delta^2_{T_{1,0}} = && Z_{4 5} \cdot Z^2_{1 2} \cdot Z^4_{3 5}
\cdot (Z^2_{1 3})^{Z^2_{3 5}Z^2_{1 2}} \cdot
(Z_{2 3})^{Z^2_{3 5}} \cdot (Z_{2 3})^{Z^{-2}_{1 2} Z^2_{3 4}} \cdot \\
&& \cdot (Z^4_{1 5})^{Z^2_{1 2}} \cdot Z^4_{3 4} \cdot (Z_{4
5})^{Z^{-2}_{1 4} Z^2_{1 2}}.
\end{eqnarray*}
\end{lem}

\begin{lem}\label{t20}
Let $T_{2,0}$ be a conic-line arrangement composed of two  conics which are  tangent at
 two points, and two lines which are tangent to the same conic (see Figure \ref{fig_t20}).

\begin{figure}[!ht]
\epsfysize=5.5cm
\centerline{\epsfbox{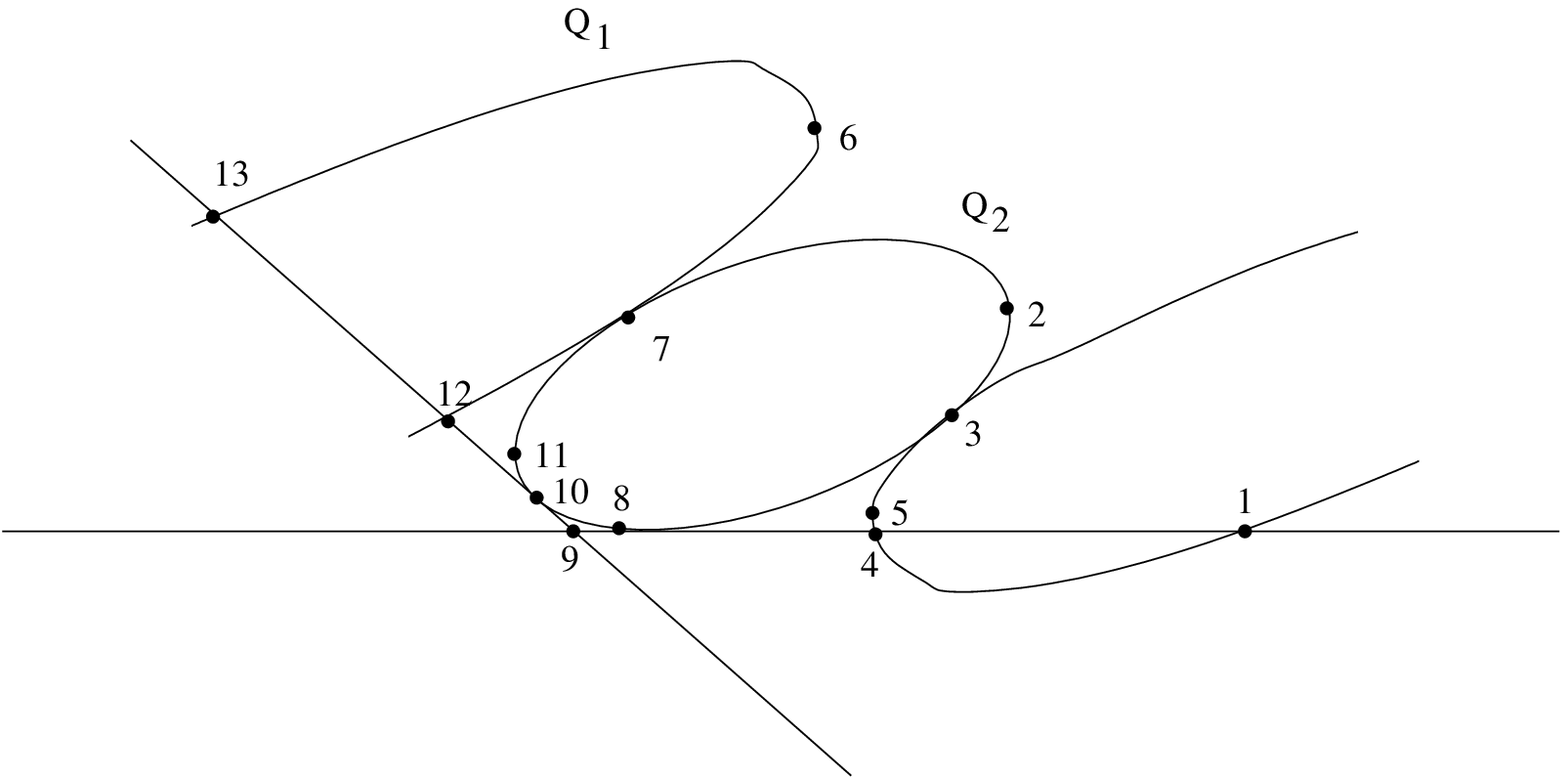}}
\caption{The arrangement $T_{2,0}$}\label{fig_t20}
\end{figure}

Then its BMF is:
\begin{eqnarray*}
\Delta^2_{T_{2,0}} = & & Z^2_{2 3} \cdot Z_{5 6} \cdot Z^4_{4 6} \cdot Z^2
_{2 3} \cdot
(Z_{3 4})^{Z^{2}_{2 3}Z^2_{4 6}} \cdot (Z_{3 4})^{Z^{-2}_{2 3} Z^2_{4 5}} \cdot Z^4_{4 5} \cdot
(Z^4_{2 6})^{Z^2_{2 3}} \cdot \\
& & (Z^2_{1 2})^{Z^2_{2 6}Z^2_{2 3}} \cdot Z^4_{1 6} \cdot  (Z_{5 6})^{Z^{-2}_{1 5}Z^{-2}_{2 5}Z^{2}_{2 3}}
\cdot (Z^2_{1 4})^{Z^{2}_{1 2}Z^{2}_{2 3}} Z^2_{1 3}.
\end{eqnarray*}
\end{lem}

Now, let us consider the conic-line arrangement $T_{n,0}$,
which is composed of two tangent conics (which are tangent
to each other at two points) and $n$ lines which are tangent to the same conic.

\begin{figure}[!ht]
\epsfysize=7cm
\centerline{\epsfbox{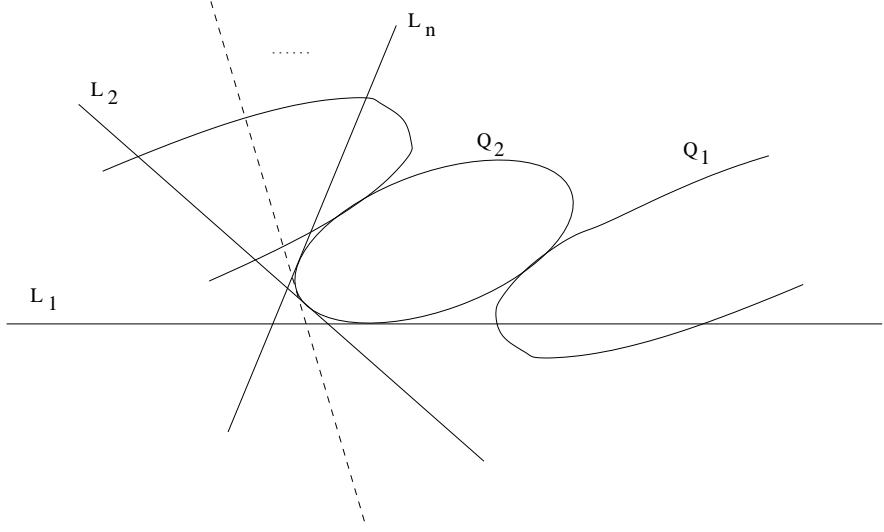}}
\caption{The arrangement $T_{n,0}$}\label{fig_tn0}
\end{figure}

\begin{prop}\label{tn0}
Let $Q_1,Q_2$ be two tangent conics in $\cpt$ (which are tangent 
to each other at two points) and let $L_1,\dots,
L_n$ be $n$ lines tangent to $Q_2$ (see Figure \ref{fig_tn0}).
Denote
$$T_{n,0}=Q_1\cup Q_2\cup \bigcup_{i=1}^n L_i.$$

Then its BMF is:

\begin{eqnarray*}
\Delta^2_{T_{n,0}}& =& Z_{n+2,n+3} \cdot \left(Z_{n,n+1}\right)^{Z^2_{n-1,
n} Z^2_{n+1,n+3}} \cdot \left(Z_{n,n+1} \right)^{Z^2_{n-1,n} Z^{-2}_{n+1,n+2}} \cdot \\
& & \cdot \tilde{Z}_{n+2,n+3} \cdot Z^4_{n+1,n+3} \cdot Z^4_{n+1,n+2} \cdot
\left(Z^4_{n-1,n+3}\right)^{Z^2_{n-1,n}} \cdot
Z^4_{n+2,n+4} \cdot \\
& & \cdot \prodl^{n-2}_{i=1} Z^4_{i,n+3}  \cdot
\prodl^{n-2}_{i=1} (Z^2_{i,n-1})^{Z^2_{n-1,n+3} Z^2_{n-1,n}} \cdot
\prodl^{n-2}_{i=1} \tilde{Z}^2_{i,n+4} \cdot \\
& & \cdot \prodl_{1\leq i < j\leq n-2} (Z^2_{i,j})^{Z^2_{j,n+3}}
\cdot \prodl^{n-2}_{i=1} \tilde{Z}^2_{i,n} \cdot
\prodl^{n-2}_{i=1} \tilde{Z}
^2_{i,n+1} \cdot  Z^2_{n-1,n} \cdot \\
& & \cdot \left(\bar{Z}^2_{n,n+4}\right)^{Z^2_{n-1,n}} \cdot \tilde{Z}^2_{n+1,n+4} \cdot
\left(Z^2_{n-1,n+4}\right)^{Z^2_{n-1,n} Z^{-2}_{n+3,n+4}}
\end{eqnarray*}

The skeletons of the braids $\tilde{Z}_{n+2,n+3}$,
$\tilde{Z}^2_{i,n+4}$,  $\tilde{ Z}^2_{i,n}$,
$\tilde{Z}^2_{i,n+1}$ and $\tilde{Z}^2_{n+1,n+4}$ appear in Figure
\ref{fig_bmf_tn0}.

\begin{figure}[!ht]
\epsfysize=9cm \centerline{\epsfbox{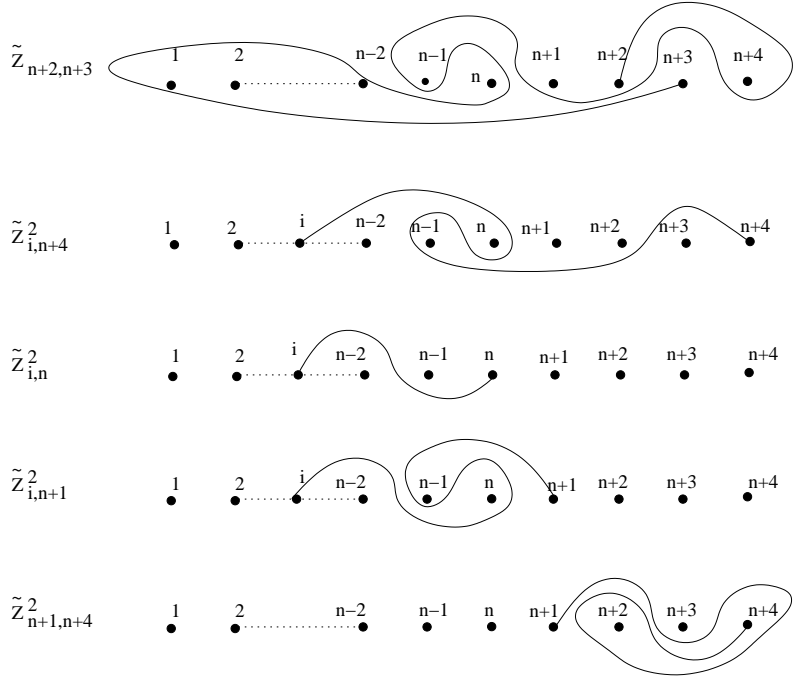}}
\caption{Skeletons for some braids appearing in  the BMF of
$T_{n,0}$}\label{fig_bmf_tn0}
\end{figure}

\end{prop}

Before the proof, we have the following remark.
\begin{rem}
The order of the factors in a BMF $\Delta^2$ is fixed according to
the order of the locations of the singularities in the curve (see
\cite{MoTe4}). Since our goal is to compute fundamental groups,
the order of the factors is irrelevant. Hence, we  list here the
monodromy factors without preserving their original order (we join
together sets of monodromies of the same type), and concentrate on
finding the relations in the group by applying the van Kampen
Theorem on the monodromies (a similar convention was also applied
in \cite{AmOg}, since we were also interested in the presentation
and not in the actual order in the factorization).
\end{rem}

\begin{proof}[Proof of Proposition \ref{tn0}]

The first four factors in $\Delta^2_{T_{n,0}}$ correspond to the four branch
 points in the arrangement.

The next five factors correspond to the tangent points in the arrangement
(two of them correspond to the tangency points between the two conics, and
the others correspond to the tangency points between one of
the conics and the lines).

The last nine factors correspond to the intersections of
the lines with the second conic, and between the lines themselves.
\end{proof}

\subsection{The corresponding fundamental groups}

In this section, we compute the projective fundamental groups
$\pi_1(\cpt-T_{n,0})$, using the BMFs we computed in the previous
section (as in \cite{AmGaTe}). The fundamental groups of the
complements of $T_{0,0}$, $T_{1,0}$ and $T_{2,0}$ have been
computed in \cite{AmGaTe} (Propositions 1.1, 1.2, 1.3), therefore
we quote here the results.

\begin{prop}
Let $T_{0,0}$ be the curve defined in Lemma \ref{t00} (see Figure \ref{fig_t00}). Then:
$$\pi_1(\cpt-T_{0,0}) \cong \langle x_1,x_2\ |\ (x_1 x_2)^2=(x_2 x_1)^2=e \rangle \cong \Z * \Z/2.$$
\end{prop}

\begin{prop}
Let $T_{1,0}$ be the curve defined in Lemma \ref{t10} (see Figure
\ref{fig_t10}). Then:
$$\pi_1(\cpt-T_{1,0}) \cong \langle x_1,x_2\ |\ (x_1 x_2)^2=(x_2 x_1)^2 \rangle.$$
\end{prop}

\begin{prop}
Let $T_{2,0}$ be the curve defined in Lemma \ref{t20} (see
Figure \ref{fig_t20}). Then:
$$\pi_1(\cpt-T_{2,0}) \cong \left\langle
\begin{array}{c|c}
x_1,x_2,x_3 & (x_2 x_3)^2=(x_3 x_2)^2, (x_1 x_3)^2=(x_3 x_1)^2,\\
            & [x_1, x_2]=[x_2, x_3 x_1 x_3^{-1}]=e
\end{array}
\right\rangle.$$
\end{prop}

\bigskip

Now we proceed to the general case:
\begin{thm}\label{presentation_Tn0}
Let $Q_1,Q_2$ be two tangent conics in $\cpt$ (which are tangent
to each other at two points) and let $L_1,\dots, L_n$ be $n$ lines tangent to
$Q_2$ (see Figure \ref{fig_tn0}). Denote
$$T_{n,0}=Q_1\cup Q_2\cup \bigcup_{i=1}^n L_i,$$
as in Proposition \ref{tn0}. Then:
\begin{tiny}
$$\pi_1(\C\P^2-T_{n,0}) \cong \left\langle \begin{array}{c} \ \\ \  \\
x_1,\dots,x_n,x_{n+2}
\\ \  \\ \  \end{array} \left| \begin{array}{cl}
(x_n x_{n+2})^2= (x_{n+2} x_n)^2 & \\
(x_i x_{n+2})^2= (x_{n+2} x_i)^2, &   1 \leq i \leq n-1 \\
\ [x_{i}, x_n]=[x_i, x_{n+2} x_n x_{n+2}^{-1}]=e, &  1 \leq i\leq n-1 \\
\ [x_i, x_{n+2} x_j x_{n+2}^{-1}]=e, & 1 \leq i < j  \leq n-1
\end{array} \right.
\right\rangle .
$$
\end{tiny}

The generators $x_1, \dots, x_{n-1}$ correspond to the lines,
and the generators $x_n$ and $x_{n+2}$ correspond to the two conics.
\end{thm}

\begin{proof}
Applying the van Kampen Theorem on $\Delta_{T_{n,0}}^2$ (computed
in Proposition
\ref{tn0}), we get the following presentation for $\pi_1(\C\P^2-T_{n,0})$:

\noindent
Generators: $\{ x_1,x_2,\dots ,x_{n+4} \}$ \\
Relations:
\begin{tiny}
\begin{enumerate}
\item $x_{n+2} = x_{n+3}$

\item $x_n x_{n-1} x_n x_{n-1}^{-1} x_n^{-1} = x_{n+3} x_{n+1} x_{n+3}^{-1}$

\item $x_n x_{n-1} x_n x_{n-1}^{-1} x_n^{-1} = x_{n+1}^{-1} x_{n+2}^{-1} x
_{n+1} x_{n+2} x_{n+1}$

\item $x_{n+3} = x_{1}^{-1} x_2^{-1} \cdots x_{n-2}^{-1} x_n x_{n-1}^{-1}
x_{n}^{-1} x_{n+3}^{-1}
x_{n+4} x_{n+3} x_{n+2} x_{n+3}^{-1} x_{n+4}^{-1} x_{n+3} x_{n} x_{n-1} x_{n}^{-1} x_{n-2}
\cdots x_{2} x_{1}$

\item $(x_{n+1} x_{n+3})^2 = (x_{n+3} x_{n+1})^2$

\item $(x_{n+1} x_{n+2})^2 = (x_{n+2} x_{n+1})^2$

\item $(x_n x_{n-1} x_n^{-1} x_{n+3})^2 = (x_{n+3} x_{n} x_{n-1} x_{n}^{-1})^2$

\item $(x_{n+2} x_{n+4})^2 = (x_{n+4} x_{n+2})^2$

\item $(x_{i} x_{n+3})^2 = (x_{n+3} x_{i})^2$, where $1 \leq i \leq n-2$

\item $[x_{i}, x_{n+3} x_n x_{n-1} x_{n}^{-1} x_{n+3}^{-1}]=e$, where $1 \leq i\leq n-2$

\item $[x_{n+3}^{-1} x_{n+4} x_{n+3}, x_{n-1}^{-1} x_{n} x_{n-1} x_{n-2} \cdots x_{i+1}
x_{i} x_{i+1}^{-1} \cdots x_{n-2}^{-1} x_{n-1}^{-1} x_{n}^{-1} x_{n-1}]=e$,
where $1 \leq i\leq n-2$

\item $[x_i, x_{n+3} x_j x_{n+3}^{-1}]=e$, where $1 \leq i < j  \leq n-2$

\item $[x_n, x_{n-2} \cdots x_{i+1} x_i x_{i+1}^{-1} \cdots x_{n-2}^{-1}]=
e $ where $1 \leq i\leq n-2$.

\item $[x_n x_{n-1}^{-1} x_{n}^{-1} x_{n+1} x_{n} x_{n-1} x_{n}^{-1},
x_{n-2}  \cdots x_{i+1} x_{i} x_{i+1}^{-1} \cdots  x_{n-2}^{-1}]=e$, where
 $1 \leq i\leq n-2$

\item $[x_{n-1}, x_n]=e$

\item $[x_n x_{n-1} x_{n} x_{n-1}^{-1} x_{n}^{-1}, x_{n+1}^{-1} x_{n+2}^{-1}
 x_{n+3}^{-1}
x_{n+4} x_{n+3} x_{n+2} x_{n+1}]=e$

\item $[x_{n+1}, x_{n+2}^{-1} x_{n+4}^{-1} x_{n+2}^{-1} x_{n+4} x_{n+2} x_{n
+4} x_{n+2}]=e$

\item $[x_{n} x_{n-1} x_{n}^{-1}, x_{n+3}^{-1} x_{n+4} x_{n+3}]=e$

\item $x_{n+4} x_{n+3} x_{n+2} x_{n+1} x_n \cdots x_1=e$ \mbox{(Projective relation)}.
\end{enumerate}
\end{tiny}

By Relation (1), the generator $x_{n+3}$ can be omitted, and we
replace it with $x_{n+2}$ in any place it appears.
Moreover, using Relation (15), we can simplify many relations, and we get the following presentation. 

\noindent
Generators: $\{ x_1,\dots, x_{n+2}, x_{n+4} \}$ \\
Relations:
\begin{tiny}
\begin{enumerate}
\item $x_n = x_{n+2} x_{n+1} x_{n+2}^{-1}$

\item $x_n = x_{n+1}^{-1} x_{n+2}^{-1} x_{n+1} x_{n+2} x_{n+1}$

\item $x_{n+2} = x_{1}^{-1} x_2^{-1} \cdots x_{n-2}^{-1} x_{n-1}^{-1} x_{n
+2}^{-1} x_{n+4} x_{n+2} x_{n+4}^{-1} x_{n+2} x_{n-1} x_{n-2} \cdots x_{2} x_{1}$

\item $(x_{n+1} x_{n+2})^2 = (x_{n+2} x_{n+1})^2$

\item $(x_{n-1} x_{n+2})^2 = (x_{n+2} x_{n-1})^2$

\item $(x_{n+2} x_{n+4})^2 = (x_{n+4} x_{n+2})^2$

\item $(x_{i} x_{n+2})^2 = (x_{n+2} x_{i})^2$, where $1 \leq i \leq n-2$

\item $[x_{n-1}, x_{n}]=e$

\item $[x_{i}, x_{n+2} x_{n-1} x_{n+2}^{-1}]=e$, where $1 \leq i\leq n-2$

\item $[x_{n+2}^{-1} x_{n+4} x_{n+2},x_{n} x_{n-2} \cdots x_{i+1} x_{i}
x_{i+1}^{-1} \cdots x_{n-2}^{-1} x_{n}^{-1}]=e$, where $1 \leq i\leq n-2$

\item $[x_i, x_{n+2} x_j x_{n+2}^{-1}]=e$, where $1 \leq i < j  \leq n-2$

\item $[x_n, x_{n-2} \cdots x_{i+1} x_i x_{i+1}^{-1} \cdots x_{n-2}^{-1}]=
e$, where $1 \leq i\leq n-2$

\item $[x_{n-1}^{-1} x_{n+1} x_{n-1},
x_{n-2}  \cdots x_{i+1} x_{i} x_{i+1}^{-1} \cdots  x_{n-2}^{-1}]=e$, where
 $1 \leq i\leq n-2$

\item $[x_n, x_{n+1}^{-1} x_{n+2}^{-2} x_{n+4} x_{n+2}^2 x_{n+1}]=e$

\item $[x_{n+1}, x_{n+2}^{-1} x_{n+4}^{-1} x_{n+2}^{-1} x_{n+4} x_{n+2} x_{n
+4} x_{n+2}]=e$

\item $[x_{n-1}, x_{n+2}^{-1} x_{n+4} x_{n+2}]=e$

\item $x_{n+4} x_{n+2}^2 x_{n+1} x_n \cdots x_1=e$.
\end{enumerate}
\end{tiny}

By Relation (6), Relation (15) is simplified to $[x_{n+1},
x_{n+4}]=e$. We combine Relations (4)--(7) together, and Relations
(9) and (11) together. So, we get the following set of relations:

\begin{tiny}
\begin{enumerate}
\item $x_n = x_{n+2} x_{n+1} x_{n+2}^{-1}$

\item $x_n = x_{n+1}^{-1} x_{n+2}^{-1} x_{n+1} x_{n+2} x_{n+1}$

\item $x_{n+2} = x_{1}^{-1} x_2^{-1} \cdots x_{n-2}^{-1} x_{n-1}^{-1} x_{n
+2}^{-1}
x_{n+4} x_{n+2} x_{n+4}^{-1} x_{n+2} x_{n-1} x_{n-2} \cdots x_{2} x_{1}$

\item $(x_{i} x_{n+2})^2 = (x_{n+2} x_{i})^2$, where $i=1, \dots, n-1, n
+1, n+4$

\item $[x_{n-1}, x_{n}]=e$

\item $[x_{n+2}^{-1} x_{n+4} x_{n+2},x_{n} x_{n-2} \cdots x_{i+1} x_{i}
x_{i+1}^{-1} \cdots x_{n-2}^{-1} x_{n}^{-1}]=e$, where $1 \leq i\leq n-2$

\item $[x_i, x_{n+2} x_j x_{n+2}^{-1}]=e$, where $1 \leq i < j  \leq n-1$

\item $[x_n, x_{n-2} \cdots x_{i+1} x_i x_{i+1}^{-1} \cdots x_{n-2}^{-1}]=
e $, where $1 \leq i\leq n-2$

\item $[x_{n-1}^{-1} x_{n+1} x_{n-1}, x_{n-2}  \cdots x_{i+1} x_{i} x_{i+1}^
{-1} \cdots  x_{n-2}^{-1}]=e$,
where $1 \leq i\leq n-2$

\item $[x_n, x_{n+1}^{-1} x_{n+2}^{-2} x_{n+4} x_{n+2}^2 x_{n+1}]=e$

\item $[x_{n+1}, x_{n+4}]=e$

\item $[x_{n-1}, x_{n+2}^{-1} x_{n+4} x_{n+2}]=e$

\item $x_{n+4} x_{n+2}^2 x_{n+1} x_n \cdots x_1=e$.
\end{enumerate}
\end{tiny}

We treat Relations (8). For  $i=n-2$, we get $[x_n, x_{n-2}]=e$.
Then, for   $i=n-3$, we get $[x_n, x_{n-2} x_{n-3}
x_{n-2}^{-1}]=e$. By the relation $[x_n,x_{n-2}]=e$, we can simplify it to
$[x_n, x_{n-3}]=e$. If we continue in a decreasing order of the
indices  $i=n-4, \dots, 1$, we get: $[x_n, x_{i}]=e$ for $1 \leq i
\leq n-2$. Now we can combine this set with Relation (5).

In a similar way, we can simplify Relations (9) to get for $1 \leq
i\leq n-2$:
$$[x_{n-1}^{-1} x_{n+1} x_{n-1}, x_{i}]=e.$$
This simplification gives us the following list of relations:

\begin{tiny}
\begin{enumerate}
\item $x_n  = x_{n+2} x_{n+1} x_{n+2}^{-1}$

\item $x_n = x_{n+1}^{-1} x_{n+2}^{-1} x_{n+1} x_{n+2} x_{n+1}$

\item $x_{n+2} = x_{1}^{-1} x_2^{-1} \cdots x_{n-2}^{-1} x_{n-1}^{-1}
x_{n+2}^{-1} x_{n+4} x_{n+2} x_{n+4}^{-1} x_{n+2}
x_{n-1} x_{n-2} \cdots x_{2} x_{1}$

\item $(x_{i} x_{n+2})^2 = (x_{n+2} x_{i})^2$, where $i=1, \dots, n-1, n
+1, n+4$

\item $[x_{i}, x_{n}]=e$, where  $1 \leq i\leq n-1$

\item $[x_{n+2}^{-1} x_{n+4} x_{n+2},x_{n} x_{n-2} \cdots x_{i+1} x_{i}
x_{i+1}^{-1} \cdots x_{n-2}^{-1} x_{n}^{-1}]=e$, where $1 \leq i\leq n-2$

\item $[x_i, x_{n+2} x_j x_{n+2}^{-1}]=e$, where $1 \leq i < j  \leq n-1$

\item $[x_{n-1}^{-1} x_{n+1} x_{n-1}, x_{i}]=e$, where $1 \leq i\leq n-2$

\item $[x_n, x_{n+1}^{-1} x_{n+2}^{-2} x_{n+4} x_{n+2}^2 x_{n+1}]=e$

\item $[x_{n+1}, x_{n+4}]=e$

\item $[x_{n-1}, x_{n+2}^{-1} x_{n+4} x_{n+2}]=e$

\item $x_{n+4} x_{n+2}^2 x_{n+1} x_n \cdots x_1=e$.
\end{enumerate}
\end{tiny}

By substituting Relation (1) in Relations (4) for $i=n+1$, we get
$(x_{n} x_{n+2})^2= (x_{n+2} x_{n})^2$. Since the relation
$(x_{n+1} x_{n+2})^2= (x_{n+2} x_{n+1})^2$ is known in Relations
(4), Relation (2) is redundant (since comparing Relations (1) and
(2) yields a known relation). Now we substitute Relation (1) in
Relation (9) and we get: 

{\begin{tiny}
\begin{eqnarray*}
& e=[x_{n+1} x_n x_{n+1}^{-1}, x_{n+2}^{-2} x_{n+4} x_{n+2}^2]=
[(x_{n+2}^{-1} x_n x_{n+2}) x_n (x_{n+2}^{-1} x_n^{-1} x_{n+2}), x_{n+2}^{-2} x_{n+4} x_{n+2}^2]=\\
& = [ x_n x_{n+2} x_n x_{n+2}^{-1} x_n^{-1}, x_{n+2}^{-1} x_{n+4} x_{n+2}]=[x_{n+2}^{-1} x_n x_{n+2}, x_{n+2}^{-1} x_{n+4} x_{n+2}]=[x_n,x_{n+4}].
\end{eqnarray*}
\end{tiny}}

By Relation (12), we substitute $x_{n+2}^2 x_{n+1} x_n \cdots x_1=x^{-1}_{n+4}$ in $[x_n,x_{n+4}]=e$.
By Relations (1) and (5), this relation can be simplified to $(x_{n+2} x_{n})^2= (x_{n} x_{n+2})^2$,
which is known, and hence Relation (9) is redundant. So we get the following set of relations:

\begin{tiny}
\begin{enumerate}
\item $x_n  = x_{n+2} x_{n+1} x_{n+2}^{-1}$

\item $x_{n+2} = x_{1}^{-1} x_2^{-1} \cdots x_{n-2}^{-1} x_{n-1}^{-1}
x_{n+2}^{-1} x_{n+4} x_{n+2} x_{n+4}^{-1} x_{n+2}
x_{n-1} x_{n-2} \cdots x_{2} x_{1}$

\item $(x_{i} x_{n+2})^2 = (x_{n+2} x_{i})^2$, where $i=1, \dots, n+1, n
+4$

\item $[x_{i}, x_{n}]=e$, where  $1 \leq i\leq n-1$

\item $[x_{n+2}^{-1} x_{n+4} x_{n+2},x_{n} x_{n-2} \cdots x_{i+1} x_{i}
x_{i+1}^{-1} \cdots x_{n-2}^{-1} x_{n}^{-1}]=e$, where $1 \leq i\leq n-2$

\item $[x_i, x_{n+2} x_j x_{n+2}^{-1}]=e$, where $1 \leq i < j  \leq n-1$

\item $[x_{n-1}^{-1} x_{n+1} x_{n-1}, x_{i}]=e$, where $1 \leq i\leq n-2$

\item $[x_{n+1}, x_{n+4}]=e$

\item $[x_{n-1}, x_{n+2}^{-1} x_{n+4} x_{n+2}]=e$

\item $x_{n+4} (x_{n+2} x_{n})^2 x_{n-1} \cdots x_1=e$.
\end{enumerate}
\end{tiny}

We start by simplifying  Relations (5)  for $i=n-2$: 
\begin{eqnarray*}
e=[x_{n+2}^{-1} x_{n+4} x_{n+2}, x_n x_{n-2} x_{n}^{-1}]
\stackrel{ (4),i=n-2}{=}
[x_{n+4}, x_{n+2} x_{n-2} x_{n+2}^{-1}].
\end{eqnarray*}
This relation can be rewritten, by using Relation (10), as:
\begin{eqnarray*}
x_n x_{n+2} x_{n} x_{n-1} \cdots x_{1}  x_{n+2} x_{n-2} x_{n+2}^{-1} =
x_{n-2} x_n x_{n+2} x_{n} x_{n-1} \cdots x_{1}.
\end{eqnarray*}
Now, by Relations (3) (for $i=n-2$), we get:
\begin{eqnarray*}
x_n x_{n+2} x_{n} x_{n-1} x_{n-2} x_{n+2} x_{n-2} x_{n+2}^{-1} x_{n-2}^{-1}=
x_{n-2} x_n x_{n+2} x_{n} x_{n-1}.
\end{eqnarray*}
By Relations (6) (for $i=n-2$), we get:
\begin{eqnarray*}
x_{n+2} x_{n} x_{n+2}^{-1} x_{n-2}  =  x_{n-2} x_{n+2} x_{n} x_{n+2}^{-1}.
\end{eqnarray*}
In a similar way, we get  $[x_{i},x_{n+2} x_{n} x_{n+2}^{-1}]=e$,
for $1 \leq i \leq n-3$. As done above, we substitute this time
Relation (10) in Relation (9) to get: $[x_{n-1},x_{n+2} x_{n}
x_{n+2}^{-1}]=e$. We combine both induced relations in Relations
 (5) below, to get:

\begin{tiny}
\begin{enumerate}
\item $x_n  = x_{n+2} x_{n+1} x_{n+2}^{-1}$

\item $x_{n+2} = x_{1}^{-1} x_2^{-1} \cdots x_{n-2}^{-1} x_{n-1}^{-1}
x_{n+2}^{-1} x_{n+4} x_{n+2} x_{n+4}^{-1} x_{n+2}
x_{n-1} x_{n-2} \cdots x_{2} x_{1}$

\item $(x_{i} x_{n+2})^2 = (x_{n+2} x_{i})^2$, where $i=1, \dots, n+1, n+4$

\item $[x_{i}, x_{n}]=e$, where  $1 \leq i\leq n-1$

\item $[x_{i},x_{n+2} x_{n} x_{n+2}^{-1}]=e$, where $1 \leq i\leq n-1$

\item $[x_i, x_{n+2} x_j x_{n+2}^{-1}]=e$, where $1 \leq i < j  \leq n-1$

\item $[x_{n-1}^{-1} x_{n+1} x_{n-1}, x_{i}]=e$, where $1 \leq i\leq n-2$

\item $[x_{n+1}, x_{n+4}]=e$

\item $x_{n+4} (x_{n+2} x_{n})^2 x_{n-1} \cdots x_1=e$.
\end{enumerate}
\end{tiny}

We simplify now Relation (8) by substituting Relation (1) and
Relation (9): $$x_{n+2}^{-1} x_n x_{n+2} \cdot  (x_{n+2} x_{n})^2
x_{n-1} \cdots x_1=  (x_{n+2} x_{n})^2 x_{n-1} \cdots x_1 \cdot
x_{n+2}^{-1} x_n x_{n+2}.$$

Using Relations (3) for $i=n$, we get $[x_{n+2} x_n x_{n+2}^{-1},
x_{n-1} \cdots x_1]=e$, which is known already by Relations (5),
and therefore is redundant. Moreover, using Relations (1), (3) and
(4), Relations (7) (for $i=n-2$) can be translated to:
\begin{eqnarray*}
&&e=[x_{n+1}, x_{n-1} x_{n-2} x_{n-1}^{-1}]=[ x_{n+2}^{-1} x_{n} x_{n+2},
x_{n-1} x_{n-2} x_{n-1}^{-1}]=\\
&& =[ x_n x_{n+2} x_{n} x_{n+2}^{-1} x_n^{-1}, x_{n-1} x_{n-2} x_{n-1}^{-1
}]=[x_{n+2} x_{n} x_{n+2}^{-1},x_{n-1} x_{n-2} x_{n-1}^{-1}].
\end{eqnarray*}
This relation is known already by Relations (5). If we check the
other cases $i=n-3, \dots, 1$, we get $[x_{n+2} x_{n}
x_{n+2}^{-1},x_{n-1} x_{i} x_{n-1}^{-1}]$, for $1 \leq i \leq
n-3$. These relations are known by Relations (5). Therefore,
Relations (7) are redundant.

Notice that it is possible to combine Relations (5) and (6) to the general form
$[x_{i}, x_{n+2} x_{j} x_{n+2}^{-1}]=e$, where $1 \leq i < j \leq n$. So we get:

\begin{enumerate}
\item $x_n  = x_{n+2} x_{n+1} x_{n+2}^{-1}$

\item $x_{n+2} = x_{1}^{-1} x_2^{-1} \cdots x_{n-2}^{-1} x_{n-1}^{-1}
x_{n+2}^{-1} x_{n+4} x_{n+2} x_{n+4}^{-1} x_{n+2}
x_{n-1} x_{n-2} \cdots x_{2} x_{1}$

\item $(x_{i} x_{n+2})^2 = (x_{n+2} x_{i})^2$, where $i=1, \dots, n+1, n
+4$

\item $[x_{i}, x_{n}]=e$, where  $1 \leq i\leq n-1$

\item $[x_i, x_{n+2} x_j x_{n+2}^{-1}]=e$, where $1 \leq i < j  \leq n$

\item $x_{n+4} (x_{n+2} x_{n})^2 x_{n-1} \cdots x_1=e$.
\end{enumerate}

Relation  $(x_{n+1} x_{n+2})^2 = (x_{n+2} x_{n+1})^2$ (from
Relations (3)) can be  translated to $(x_{n} x_{n+2})^2 = (x_{n+2}
x_{n})^2$ by substituting Relation (1), and therefore is
redundant. 

Our goal now is to omit $x_{n+4}$. By Relations (3) for
$i=n$ and Relations (4) and (6), Relation (3) for $i=n+4$ can be
rewritten as (by substituting $x_{n+4}$ from Relation (2))
$$(x_{n+2} \cdot x_{n} x_{n-1} \cdots x_{1})^2=(x_n x_{n-1} \cdots
x_{1} \cdot x_{n+2})^2.$$ Using Relation (6) again, Relation (2)
becomes:
\begin{eqnarray*}
 x_{n+2} = && x_{1}^{-1} \cdots  x_{n-1}^{-1} x_{n+2}^{-1}
(x_{1}^{-1} \cdots x_{n-1}^{-1}
(x_{n+2} x_n)^{-2}) x_{n+2} \cdot \\
&& \cdot ((x_{n+2} x_{n})^2 x_{n-1} \cdots x_1) x_{n+2} x_{n-1}
\cdots x_{1},
\end{eqnarray*}
which can be written as $(x_{n+2} \cdot x_{n-1} \cdots x_{1})^2=(x_{n-1}
\cdots x_{1} \cdot x_{n+2})^2$.
By Relation (1), the generator $x_{n+1}$ is omitted. By Relation (6), $x_{n+4}$  can be omitted too.
Therefore, we have the following presentation:\\
Generators: $\{ x_1,\dots, x_{n}, x_{n+2} \}$ \\
Relations:
\begin{enumerate}
\item $(x_{n+2} \cdot x_{n-1} \cdots x_{1})^2=(x_{n-1} \cdots
x_{1} \cdot x_{n+2})^2$

\item $(x_{i} x_{n+2})^2 = (x_{n+2} x_{i})^2$, where $1 \leq i
\leq n$

\item $(x_{n+2} \cdot x_{n} x_{n-1} \cdots x_{1})^2=(x_n x_{n-1}
\cdots x_ {1} \cdot x_{n+2})^2$

\item $[x_{i}, x_{n}]=e$, where  $1 \leq i\leq n-1$

\item $[x_i, x_{n+2} x_j x_{n+2}^{-1}]=e$, where $1 \leq i < j  \leq n$
\end{enumerate}

Now, we show that Relations (1) and (3) are redundant. Let us
denote $b:=x_n$, $a:=x_{n-1} \cdots x_{1}$ and $x:=x_{n+2} $. In
these notations, we have the relations $[a,b]=e$ (by Relations
(4)), $(bx)^2 = (xb)^2$ (by Relations (2)), $(ax)^2 = (xa)^2$ (by
Relation (1)) and $[xbx^{-1},a]=e$ (by Relations (5)). We  simplify
Relation (3), which is $(bax)^2 = (xba)^2$. We do it as follows:

\begin{eqnarray*}
&& baxbax=xbaxba \rightarrow baxbx^{-1}xax=xbaxbx^{-1}xa  \rightarrow  \\
&& \rightarrow bxbx^{-1}axax=xbxbx^{-1}axa.
\end{eqnarray*}
The last resulting relation is redundant (that means that Relation (3) is redundant).

We show now that Relation (1) is redundant too. First, we write
this relation in detail:

\begin{eqnarray*}
&& x_{n+2} \cdot x_{n-1} x_{n-2} \cdots x_2  x_{1} \cdot x_{n+2} \cdot x_{n-1} x_{
n-2} \cdots  x_2  x_{1} =\\
&& =x_{n-1} x_{n-2} \cdots x_2  x_{1} \cdot x_{n+2} \cdot x_{n-1} x_{n-2} \cdots
 x_2 x_{1} \cdot x_{n+2}.
\end{eqnarray*}

We add the expression $x_{n+2}^{-1} x_{n+2}=e$:

\begin{tiny}
\begin{eqnarray*}
&&x_{n+2} x_{n-1} x_{n-2} \cdots x_2  x_{1} (x_{n+2} x_{n-1} x_{n+2}^{-1}) (
x_{n+2} x_{n-2}
x_{n+2}^{-1}) \cdots (x_{n+2} x_2 x_{n+2}^{-1}) x_{n+2} x_{1}=\\
&&=x_{n-1} x_{n-2} \cdots x_2  x_{1} (x_{n+2} x_{n-1} x_{n+2}^{-1}) (x_{n+
2} x_{n-2} x_{n+2}^{-1}) \cdots
(x_{n+2} x_2 x_{n+2}^{-1}) x_{n+2} x_{1} x_{n+2}.
\end{eqnarray*}
\end{tiny}

By Relations (5), $x_1$ commutes with all expressions in the brackets. Therefore
we get:

\begin{tiny}
\begin{eqnarray*}
&&x_{n+2} x_{n-1} x_{n-2} \cdots x_2 (x_{n+2} x_{n-1} x_{n+2}^{-1}) (x_{n+2}
 x_{n-2}
x_{n+2}^{-1}) \cdots (x_{n+2} x_2 x_{n+2}^{-1}) x_1 x_{n+2} x_{1}=\\
&&=x_{n-1} x_{n-2} \cdots x_2 (x_{n+2} x_{n-1} x_{n+2}^{-1}) (x_{n+2} x_{n
-2} x_{n+2}^{-1}) \cdots
(x_{n+2} x_2 x_{n+2}^{-1}) x_1 x_{n+2} x_{1} x_{n+2}.
\end{eqnarray*}
\end{tiny}

Since  by Relations (2), $(x_{1} x_{n+2})^2 = (x_{n+2}
x_{1})^2$, we get:

\begin{tiny}
\begin{eqnarray*}
&&x_{n+2} x_{n-1} x_{n-2} \cdots x_2 (x_{n+2} x_{n-1} x_{n+2}^{-1}) (x_{n+2}
 x_{n-2}
x_{n+2}^{-1}) \cdots (x_{n+2} x_3 x_{n+2}^{-1})(x_{n+2} x_2 x_{n+2}^{-1})=\\
&&=x_{n-1} x_{n-2} \cdots x_2 (x_{n+2} x_{n-1} x_{n+2}^{-1}) (x_{n+2} x_{n
-2} x_{n+2}^{-1}) \cdots (x_{n+2} x_3 x_{n+2}^{-1}) x_{n+2} x_2.
\end{eqnarray*}
\end{tiny}

Again by Relations (5),  $x_2$ commutes with all the
expressions in the brackets. Hence  we have:

\begin{tiny}
\begin{eqnarray*}
&&x_{n+2} x_{n-1} x_{n-2} \cdots x_3 (x_{n+2} x_{n-1} x_{n+2}^{-1}) (x_{n+2}
 x_{n-2} x_{n+2}^{-1}) \cdots (x_{n+2} x_3 x_{n+2}^{-1}) x_2 x_{n+2} x_2=\\
&&=x_{n-1} x_{n-2} \cdots x_3 (x_{n+2} x_{n-1} x_{n+2}^{-1}) (x_{n+2}
x_{n-2} x_{n+2}^{-1}) \cdots (x_{n+2} x_3 x_{n+2}^{-1}) x_2 x_{n+2} x_2 x_{n+2}.
\end{eqnarray*}
\end{tiny}

Since by Relations (2),  $(x_{2} x_{n+2})^2 = (x_{n+2}
x_{2})^2$, we get a
 shorter form:

\begin{tiny}
\begin{eqnarray*}
&&x_{n+2} x_{n-1} x_{n-2} \cdots x_3 (x_{n+2} x_{n-1} x_{n+2}^{-1}) (x_{n+2}
 x_{n-2} x_{n+2}^{-1})
\cdots (x_{n+2} x_{3} x_{n+2}^{-1})=\\
&&=x_{n-1} x_{n-2} \cdots x_3 (x_{n+2} x_{n-1} x_{n+2}^{-1}) (x_{n+2} x_{n
-2} x_{n+2}^{-1}) \cdots x_{n+2} x_3.
\end{eqnarray*}
\end{tiny}

We continue in the same manner for $x_3, x_4, \dots, x_{n-2}$ in order to get the
following relation:

\begin{tiny}
\begin{eqnarray*}
x_{n+2} x_{n-1}(x_{n+2} x_{n-1} x_{n+2}^{-1}) x_{n-2} x_{n+2} x_{n-2} x_{n+2
}^{-1}=
x_{n-1}(x_{n+2} x_{n-1} x_{n+2}^{-1}) x_{n-2} x_{n+2} x_{n-2}.
\end{eqnarray*}
\end{tiny}

Since $(x_{n-2} x_{n+2})^2 = (x_{n+2} x_{n-2})^2$ (by Relations
(2)), we get  $(x_ {n-1} x_{n+2})^2 = (x_{n+2} x_{n-1})^2$, which
appear already in Relations (2). Therefore, Relation  (1) is
redundant.

\medskip

We finally get the requested presentation for
$\pi_1(\C\P^2-T_{n,0})$ with $x_1, \dots, x_n, x_{n+2}$ as
generators and with the following relations:

\begin{enumerate}
\item $(x_{i} x_{n+2})^2 = (x_{n+2} x_{i})^2$, where $1 \leq i
\leq n$

\item $[x_{i}, x_{n}]=e$, where  $1 \leq i\leq n-1$

\item $[x_i, x_{n+2} x_j x_{n+2}^{-1}]=e$, where $1 \leq i < j  \leq n$,
\end{enumerate}
as needed.

\end{proof}

By the above results, we can conclude:

\begin{cor}\label{big2}
The affine and projective fundamental groups of a conic-line
arrangement, composed of two  conics, which are tangent to each
other at two points, and $n$ tangent lines, are big.
\end{cor}

\begin{proof}

We have seen that the projective fundamental group of two tangent
conics is $\langle a,b | (ab)^2=(ba)^2= e\rangle \cong \Z * \Z/2$, which is big, as we saw in Corollary \ref{big}.

Since the projective fundamental group is a quotient of the affine fundamental
group
by the projective relation, the affine fundamental group is big too.

Now, let $A$ be a conic-line arrangement with two conics
and with some additional tangent lines. Since the arrangement composed of
two tangent conics is a sub-arrangement of $A$, and the fundamental groups
of two tangent conics are big, we have that the fundamental groups of $A$ are big too.
\end{proof}

\section{The fundamental group of the complement of two tangent conics with
$n+m$ tangent lines}\label{sec4}

In this section, we compute the fundamental group of the
complement of the arrangement $T_{n,m}$ consisting of two tangent conics, $n$
 lines which are
tangent to one of the conics and $m$ lines which are  tangent to the other one.

\subsection{The braid monodromy factorizations}\label{bm_n+m}

In this subsection we present the BMF of the arrangements
$T_{1,1}$, $T_{2,1}$, $T_{2,2}$ and we shall compute the BMF of
the arrangement $T_{n,m}$ for any $n, m$. Note that we dealt with
the arrangement $T_{1,1}$ in \cite{AmGaTe}, and we give here its
BMF $\Delta^2_{T_{1,1}}$.

\begin{lem}\label{t11}
Let $T_{1,1}$ be a conic-line arrangement composed of two tangent conics (at
 two points) and
two tangent lines, each is tangent to a different conic, see Figure \ref{thm3-fig}.

\begin{figure}[h]
\epsfysize=5cm
\centerline{\epsfbox{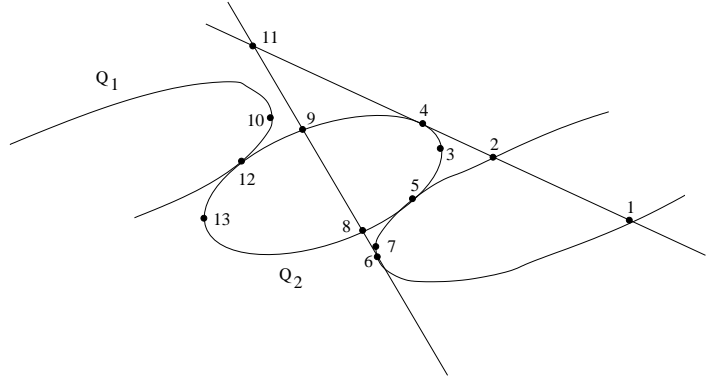}}
\caption{The arrangement $T_{1,1}$}\label{thm3-fig}
\end{figure}
Then its BMF is
\begin{eqnarray*}
\Delta^2_{T_{1,1}} = & & {Z}^2_{2 3} \cdot \bar{Z}^2_{2 4} \cdot
{(Z_{5 6})}^{Z^{-2}_{2 5}} \cdot Z_{2 5}^4 \cdot Z_{4 6}^{4} \cdot
Z_{1 3}^{ 4} \cdot (Z_{3 4})^{Z_{1 3}^2 Z_{4 6}^{2}} \cdot \\ &&
\cdot (Z^2_{1 6})^{Z^2_{1 3}} \cdot (Z_{1 5}^{2})^{Z_{1 2}^{2}}
\cdot (Z_{3 4})^{Z^2_{4 5}} \cdot Z^2_{1 2} \cdot  Z_{4 5}^{4}
\cdot Z_{5 6}.
\end{eqnarray*}
\end{lem}

In a similar way, one can compute the  BMFs of $T_{2,1}$ and
$T_{2,2}$:

\begin{lem}\label{t21}
Let $T_{2,1}$ be a conic-line arrangement composed of two tangent conics (at
 two points) and three additional  lines, two of them are tangent to the first conic,
and the third one is tangent to the second conic,
see Figure \ref{fig-t21}.

\begin{figure}[h]
\epsfysize=5cm
\centerline{\epsfbox{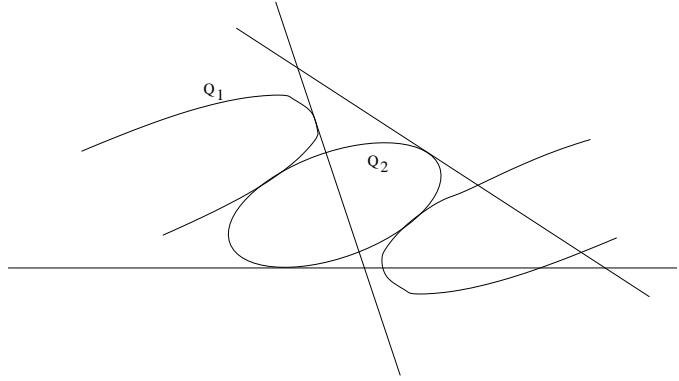}}
\caption{The arrangement $T_{2,1}$}\label{fig-t21}
\end{figure}

Then its BMF is:
\begin{eqnarray*}
\Delta^2_{T_{2,1}} & = & {Z}_{2 3} \cdot {Z}^2_{6 7} \cdot Z_{1 3}^{4}
\cdot Z_{2 4 }^{4} \cdot
Z_{5 7}^{4} \cdot (Z_{4 5})^{\bar{Z}_{2 4}^2 Z_{5 7}^{2}} \cdot (Z^2_{2 7})^
{Z^2_{2 4}Z^2_{2 3}}
\cdot (Z_{3 7}^{2})^{Z_{1 3}^{2}} \cdot \\
& & \cdot Z^4_{3 6} \cdot (Z^2_{1 7})^{Z^2_{1 3}} \cdot (Z_{4
5})^{{Z}_{3 4}^{-2} Z_{1 3}^{2}} \cdot
(Z^4_{3 4})^{Z^2_{1 3}} \cdot (Z^2_{1 5})^{Z^2_{1 3}} \cdot (Z^2_{1 5})^{Z^2_{1 3}} \cdot \\
&& \cdot (Z_{2 3})^{Z^2_{3 6} \bar{Z}_{1 3}^{-2}} \cdot {Z}^2_{4
6} \cdot (Z^2_{5 6})^{Z^2_{1 5}{Z}_{1 3}^{2}} \cdot (Z^2_{1
6})^{Z^2_{1 5}Z^2_{1 3}}.
\end{eqnarray*}
\end{lem}

\begin{lem}\label{t22}
Let $T_{2,2}$ be a conic-line arrangement composed of two tangent conics (at
 two points) and four additional lines, two of them are tangent to the first conic, and the two
others are tangent to the second conic,
see Figure \ref{fig-t22}.

\begin{figure}[h]
\epsfysize=5cm
\centerline{\epsfbox{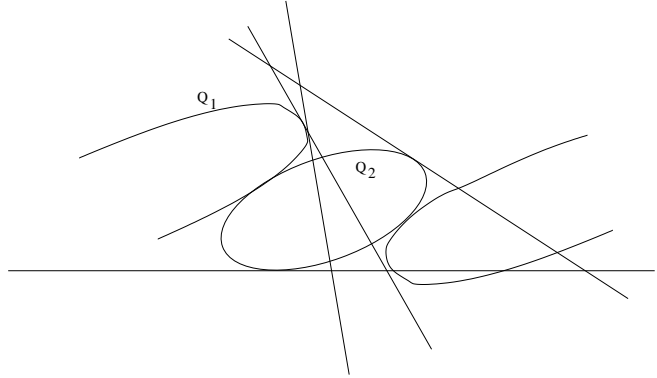}}
\caption{The arrangement $T_{2,2}$}\label{fig-t22}
\end{figure}

Then its BMF is:
\begin{eqnarray*}
\Delta^2_{T_{2,2}} & = & {Z}_{2 3} \cdot {Z}^2_{6 7} \cdot Z_{1 3}^{4}
\cdot Z_{2 4 }^{4} \cdot
Z_{5 7}^{4} \cdot (Z_{4 5})^{\bar{Z}_{2 4}^2 Z_{5 7}^{2}} \cdot (Z^2_{2 7})^
{Z^2_{2 4}Z^2_{2 3}}
\cdot (Z_{3 7}^{2})^{Z_{1 3}^{2}} \cdot \\
& & \cdot (Z^2_{6 8})^{Z^2_{6 7}} \cdot (\bar{Z}^2_{2
8})^{\bar{Z}_{6 8}^{2}} \cdot (Z^2_{3 8})^{Z^2_{3 7}Z^2_{3
5}Z^2_{1 3}} \cdot {Z^4_{2 6}} \cdot (Z^2_{1 7}
)^{{Z}_{1 3}^{2}} \cdot \\
& & \cdot (Z_{4 5})^{Z^{-2}_{5 8}{Z}_{7 8}^{-2}Z^{-2}_{3 4}Z^2_{1
3}} \cdot (Z_{5 8}^{4})^{Z_{5 7}^{2}} \cdot (Z^2_{1 8})^{Z^2_{1
7}Z^2_{1 3}} \cdot {Z} ^2_{7 8} \cdot
(Z_{1 5}^{2})^{Z_{1 3}^{2}} \cdot \\
& & \cdot (Z^4_{3 4})^{Z^2_{1 3}} \cdot (Z_{1 5}^{2})^{Z_{1
3}^{2}} \cdot (Z_{2 3})^{Z^{-2}_{1 2}\bar{Z}^2_{2 6}} \cdot
(Z^2_{4 6})^{Z^2_{4 5}} \cdot
 (Z^2_{5 6})^{Z^2_{1 5}{Z}^2_{1 3}} \cdot (Z^2_{1 6})^{Z_{1 5}^2{Z}_{1 3}^{2}}.
\end{eqnarray*}
\end{lem}

Now, we compute the BMF of $T_{n,m}$.

\begin{prop}\label{bmf-tnm}
Let $Q_1,Q_2$ be two tangent conics in $\cpt$. Consider two sets
of $n$ and $m$ lines, such that $n$ lines are tangent to $Q_2$, and
$m$  lines are tangent to $Q_1$. Denote: $T_{n,m}=Q_1\cup Q_2\cup (\bigcup_{i=1}^n L_i) \cup (\bigcup_{j=1}^m L'_j)$.
Then:
\begin{eqnarray*}
\Delta_{T_{n,m}}^2  & =  & Z_{2 3} \cdot (Z_{2 3})^{{Z}_{1 2}^{-2} \prodlim_{i=6}^{n+4} \bar{Z}^2_{2 i}} \cdot
                       {(Z_{4 5})^{\bar{Z}^2_{2 4} Z^2_{5,n+5}}} \cdot \tilde{Z}_{4 5} \cdot
                       \bar{Z}^4_{2 4} \cdot {(Z^4_{3 4})}^{Z^2_{1 3}} \cdot \\
                     & & \cdot {Z^4_{1 3}}  \cdot \prodlim_{i=6}^{n+4} \bar{Z}^4_{2 i} \cdot
                       \prodlim_{i=n+5}^{n+m+4} {(Z^4_{5 i})}^{Z^2_{5,n+5}} \cdot
                       {(Z^2_{1 5})}^{Z^2_{1 3}} \cdot \prodlim_{i=6}^{n+4} \Big{[}{(Z^2_{4 i})}^{Z^{2}_{4 5}}
                       \cdot {(Z^2_{5 i})}^{Z^2_{1 5} Z^2_{1 3}}\Big{]}  \cdot \\
                     & & \cdot {\prodlim_{i=6}^{n+4} {(Z^2_{1 i})}^{Z^2_{1 5} Z^2_{1 3}}}  \cdot
                       \prodlim_{6 \leq i < j \leq n+4}^{} {(\bar{Z}^2_{i j})}^{\bar{Z}^{2}_{2 i}}  \cdot
                       {(Z^2_{2,n+5})}^{Z^2_{2 4} Z^2_{2 3}} \cdot {(Z^2_{3,n+5})}^{Z^2_{1 3}} \cdot \\
                     & & \cdot \prodlim_{i=n+6}^{n+m+4} \Big{[}{(\bar{Z}^2_{2 i})}^{(\prodlim_{j=n+4}^{6} \bar{Z}^{-2}_{2 j})}
                       \cdot \tilde{Z}^2_{3 i} \Big{]}  \cdot {\prodlim_{i=n+6}^{n+m+4} Z^2_{n+5,i}} \cdot
                       \prodlim_{n+6 \leq i < j \leq n+m+4}^{} {(Z^2_{i j})}^{Z^{-2}_{5 i} Z^2_{5,n+5}} \cdot  \\
                     & & \cdot {\prodlim_{i=n+5}^{n+m+4} {({Z}^2_{1 i})}^{Z^{2}_{1,n+5} Z^2_{1 3}}}  \cdot
                       \prodlim_{i=6}^{n+4} \Big{[}\prodlim_{j=n+5}^{n+m+4} {({Z}^2_{i j})}^
                       {(\prodlim_{k=n+5}^{j}Z^{-2}_{k j})} \Big{]}
\end{eqnarray*}
The skeletons of  $\tilde{Z}_{4 5}$ and $\tilde{Z}^2_{3 i}$ appear
in Figure \ref{2figsnew}.

\begin{figure}[h]
\epsfysize=4.5cm \centerline{\epsfbox{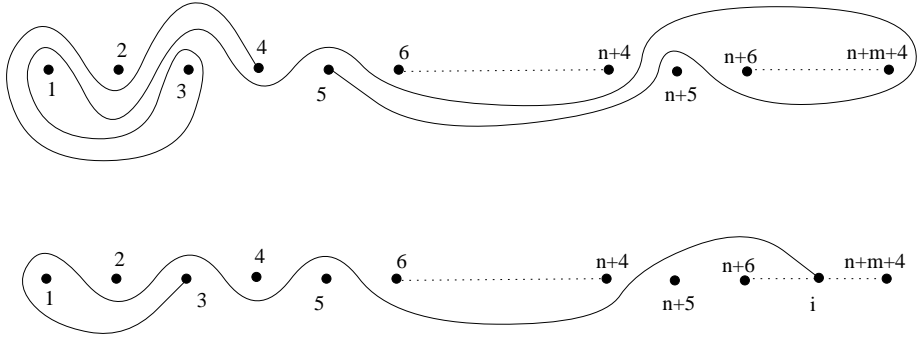}}
\caption{Skeletons for the BMF of $T_{n,m}$}\label{2figsnew}
\end{figure}

\end{prop}

\begin{proof}
By a similar technique to what we have used for computing the BMF of $T_{n,0}$, we compute the BMF of $T_{n,m}$.

We first construct the BMF of $T_{1,m}$ by generalizing the BMF of
$T_{1,1}$, $T_{1,2}$, $T_{1,3}$:
\begin{eqnarray*}
\Delta_{T_{1,m}}^2 & = & Z_{2 3} \cdot (Z_{2 3})^{{Z}_{1 2}^{-2}} \cdot (Z_{4 5})^{\bar{Z}^2_{2 4} Z^2_{5 6}} \cdot
\tilde{Z'}_{4 5} \cdot \bar{Z}^4_{2 4} \cdot {(Z^4_{3 4})}^{Z^2_{1 3}} \cdot Z^4_{1 3} \cdot
\prodlim_{i=6}^{m+5} {(Z^4_{5 i})}^{Z^2_{5 6}} \cdot \\
& &  \cdot {(Z^2_{1 5})}^{Z^2_{1 3}} \cdot {(Z^2_{2 6})}^{Z^2_{2
4} Z^2_{2 3}} \cdot {(Z^2_{3 6})}^{Z^2_{1 3}} \cdot
\Big{(}\prodlim_{i=7}^{m+5} \Big{[}\bar{Z}^2_{2 i} \cdot
\tilde{Z'}^2_{3 i} \Big{]}\Big{)} \cdot \prodlim_{i=7}^{m+5}
Z^2_{6 i} \cdot \\
& &  \cdot \prodlim_{7 \leq i < j \leq m+5}^{} {(Z^2_{i
j})}^{Z^{-2}_{5 i} Z^2_{5 6}} \cdot \prodlim_{i=6}^{m+5} {(Z^2_{1
i})}^{Z^{2}_{1 6} Z^2_{1 3}}.
\end{eqnarray*}

Then we proceed in the same manner to compute the BMF of $T_{2,m}$
and $T_{3,m}$, as presented in the following table. The table is
constructed as follows: in each row we write the factors which are
related to the same type of singularities in the different
arrangements.
$$
{\begin{tiny}
\begin{array}{|c||c|c|c|}
\hline
 & \Delta_{T_{1,m}}^2 & \Delta_{T_{2,m}}^2 & \Delta_{T_{3,m}}^2\\
\hline
\hline
(1) & Z_{2 3} &  Z_{2 3} & Z_{2 3}\\
\hline(2) & (Z_{2 3})^{{Z}_{1 2}^{-2}} & (Z_{2 3})^{{Z}_{1 2}^{-2} \bar{Z}^2_{2 6}} &
(Z_{2 3})^{{Z}_{1 2}^{-2} \bar{Z}^2_{2 6} \bar{Z}^2_{2 7}} \\
\hline(3) & (Z_{4 5})^{\bar{Z}^2_{2 4} Z^2_{5 6}} & (Z_{4 5})^{\bar{Z}^2_{2 4} Z^2_{5 7}} &
(Z_{4 5})^{\bar{Z}^2_{2 4} Z^2_{5 8}} \\
\hline(4) & {\tilde{Z}_{4 5}}' & {\tilde{Z}_{4 5}}'' & {\tilde{Z}_{4 5}}''' \\
\hline(5) & \bar{Z}^4_{2 4} \cdot {(Z^4_{3 4})}^{Z^2_{1 3}} & \bar{Z}^4_{2 4} \cdot {(Z^4_{3 4})}^{Z^2_{1 3}}
& \bar{Z}^4_{2 4} \cdot {(Z^4_{3 4})}^{Z^2_{1 3}}\\
\hline(6) & Z^4_{1 3} & Z^4_{1 3} \cdot \bar{Z}^4_{2 6} & Z^4_{1 3} \cdot \prodlim_{i=6,7}^{} {\bar{Z}^4_{2 i}}\\
\hline(7) & \prodlim_{i=6}^{m+5} {(Z^4_{5 i})}^{Z^2_{5 6}} & \prodlim_{i=7}^{m+6} {(Z^4_{5 i})}^{Z^2_{5 7}} &
\prodlim_{i=8}^{m+7} {(Z^4_{5 i})}^{Z^2_{5 8}} \\
\hline(8)& {(Z^2_{1 5})}^{Z^2_{1 3}} & {(Z^2_{1 5})}^{Z^2_{1 3}}
\cdot {(Z^2_{4 6})}^{Z^2_{4 5}} \cdot {(Z^2_{5 6})}^{Z^2_{1 5}
Z^2_{1 3}} & {(Z^2_{1 5})}^{Z^2_{1 3}} \cdot \prodlim_{i=6}^{7}
\Big{[}{(Z^2_{4 i})}^
{Z^{2}_{4 5}} \cdot {(Z^2_{5 i})}^{Z^2_{1 5} Z^2_{1 3}}\Big{]} \\
\hline(9)& - & {(Z^2_{1 6})}^{Z^2_{1 5} Z^2_{1 3}} &
\prodlim_{i=6}^{7} {(Z^2_{1 i})}^{Z^2_{1 5} Z^2_{1 3}}
\cdot {(Z^2_{6 7})}^{\bar{Z}^2_{2 6}}\\
\hline(10)& {(Z^2_{2 6})}^{Z^2_{2 4} Z^2_{2 3}} \cdot {(Z^2_{3 6})}^{Z^2_{1 3}} \cdot &
{(Z^2_{2 7})}^{Z^2_{2 4} Z^2_{2 3}} \cdot {(Z^2_{3 7})}^{Z^2_{1 3}} \cdot &
{(Z^2_{2 8})}^{Z^2_{2 4} Z^2_{2 3}} \cdot {(Z^2_{3 8})}^{Z^2_{1 3}} \cdot \\
& \cdot \prodlim_{i=7}^{m+5} \Big{[}\bar{Z}^2_{2 i} \cdot \tilde{Z_{3i}'}^2 \Big{]}&
\cdot \prodlim_{i=8}^{m+6} \Big{[}{(\bar{Z}^2_{2 i})}^{\bar{Z}^{-2}_{2 6}} \cdot \tilde{Z_{3i}''}^2 \Big{]}&
\cdot \prodlim_{i=9}^{m+7} \Big{[}{(\bar{Z}^2_{2 i})}^{\bar{Z}^{-2}_{2 7} \bar{Z}^{-2}_{2 6}}
\cdot \tilde{Z_{3i}'''}^2 \Big{]} \\
\hline(11) & \prodlim_{i=7}^{m+5} Z^2_{6 i} \cdot & \prodlim_{i=8}^{m+6} Z^2_{7 i} \cdot
& \prodlim_{i=9}^{m+7} Z^2_{8 i} \cdot \\
& \cdot \prodlim_{7 \leq i < j \leq m+5}^{} {(Z^2_{i j})}^{Z^{-2}_{5 i}
Z^2_{5 6}} & \cdot \prodlim_{8 \leq i < j \leq m+6}^{} {(Z^2_{i j})}^{Z^{-2}_{5 i} Z^2_{5 7}} &
\cdot \prodlim_{9 \leq i < j \leq m+7}^{} {(Z^2_{i j})}^{Z^{-2}_{5 i} Z^2_{5 8}} \\
\hline(12) & \prodlim_{i=6}^{m+5} {(Z^2_{1 i})}^{Z^{2}_{1 6}Z^2_{1 3}} &
\prodlim_{i=7}^{m+6} {(Z^2_{1 i})}^{Z^{2}_{1 7}Z^2_{1 3}} \cdot &
\prodlim_{i=8}^{m+7} {(Z^2_{1 i})}^{Z^{2}_{1 8} Z^2_{1 3}} \cdot \\
& & \cdot \prodlim_{j=7}^{m+6} {({Z}^2_{6 j})}^{(\prodlim_{k=7}^{j} Z^{-2}_{k j})} &
\cdot \prodlim_{i=6}^{7} \Big{[}\prodlim_{j=8}^{m+7} {({Z}^2_{i j})}^{(\prodlim_{k=8}^{j} Z^{-2}_{k j})} \Big{]} \\
\hline
\end{array}
\end{tiny}}
$$
The skeletons of  ${\tilde{Z}_{4 5}}',  {\tilde{Z}_{4 5}}'',
{\tilde{Z}_{4 5}}'''$  and $\tilde{Z'_{3i}}^2, \tilde{Z''_{3i}}^2, 
\tilde{Z'''_{3i}}^2$ appear in Figures \ref{3rows1} and
\ref{3rows2} respectively.
\begin{figure}[h]
\epsfysize=6cm \centerline{\epsfbox{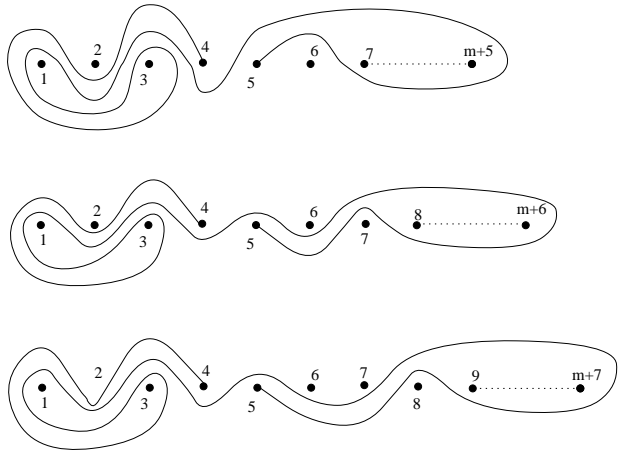}}
\caption{The skeletons of ${\tilde{Z}_{4 5}}',  {\tilde{Z}_{4
5}}'', {\tilde{Z}_{4 5}}'''$}\label{3rows1}
\end{figure}
\begin{figure}[h]
\epsfysize=6cm \centerline{\epsfbox{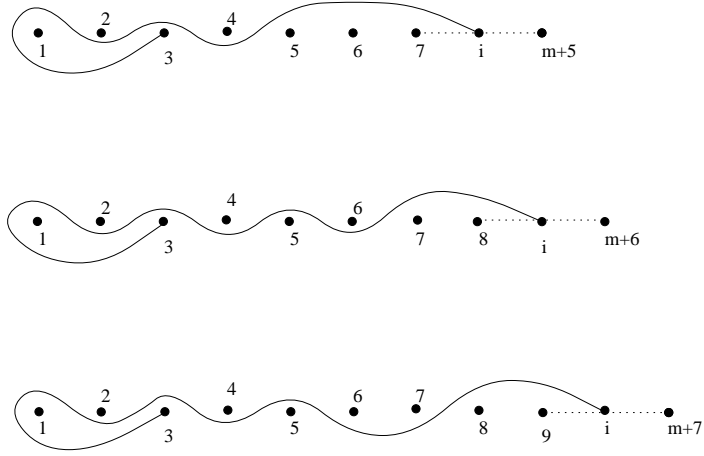}}
\caption{The skeletons of $\tilde{Z'_{3i}}^2, \tilde{Z''_{3i}}^2, 
\tilde{Z'''_{3i}}^2$}\label{3rows2}
\end{figure}

Rows (1) and (2) correspond to the branch points of $Q_1$. Rows
(3) and (4) correspond to the branch points of $Q_2$. Row (5)
corresponds to the two tangency points between $Q_1$ and $Q_2$.
Row (6) corresponds to the tangency points between $Q_1$ and the
lines which are tangent to it.  Row (7) corresponds to the
tangency points between $Q_2$ and the lines which are tangent to
it. Row (8) corresponds to the intersection points between the
lines $\{L_i\}_{i=1}^n$ and the conic $Q_2$. Row (9) corresponds
to the intersection points between the lines $\{L_i\}_{i=1}^n$
themselves. Row (10) corresponds to the intersection points
between the lines $\{L'_j\}_{j=1}^m$ with $Q_1$. Row (11)
corresponds to the intersection points between the lines
$\{L'_j\}_{j=1}^m$ themselves. Row (12) corresponds to the
intersection points between the lines $\{L_i\}_{i=1}^n$ and
$\{L'_j\}_{j=1}^m$.

In the next step, we generalize the factorization for any $n$, i.e.
we compute $\Delta^2_{T_{n,m}}$. In the following table, each row
includes the general form of the corresponding row in the previous
table:
$$
\begin{array}{|c|c|}
\hline
& \Delta_{T_{n,m}}^2 \\
\hline
\hline (1)& Z_{2 3}   \\
\hline (2)&(Z_{2 3})^{{Z}_{1 2}^{-2} \prodlim_{i=6}^{n+4} \bar{Z}^2_{2 i}}  \\
\hline (3)&(Z_{4 5})^{\bar{Z}^2_{2 4} Z^2_{5,n+5}}  \\
\hline (4)&\tilde{Z}_{4 5}  \\
\hline (5)&\bar{Z}^4_{2 4} \cdot {(Z^4_{3 4})}^{Z^2_{1 3}} \\
\hline (6)&Z^4_{1 3}  \cdot \prodlim_{i=6}^{n+4} \bar{Z}^4_{2 i}  \\
\hline (7)&\prodlim_{i=n+5}^{n+m+4} {(Z^4_{5 i})}^{Z^2_{5,n+5}}  \\
\hline (8)&{(Z^2_{1 5})}^{Z^2_{1 3}} \cdot \prodlim_{i=6}^{n+4} \Big{[}{(Z^2_{4 i})}^{Z^{2}_{4 5}}
       \cdot {(Z^2_{5 i})}^{Z^2_{1 5} Z^2_{1 3}}\Big{]}  \\
\hline (9)&\prodlim_{i=6}^{n+4} {(Z^2_{1 i})}^{Z^2_{1 5} Z^2_{1 3}}  \cdot
       \prodlim_{6 \leq i < j \leq n+4}^{} {(\bar{Z}^2_{i j})}^{\bar{Z}^{2}_{2 i}}  \\
\hline (10)&{(Z^2_{2,n+5})}^{Z^2_{2 4} Z^2_{2 3}} \cdot {(Z^2_{3,n+5})}^{Z^2_{1 3}}
       \cdot \prodlim_{i=n+6}^{n+m+4} \Big{[}{(\bar{Z}^2_{2 i})}^{(\prodlim_{j=n+4}^{6} \bar{Z}^{-2}_{2 j})} \cdot
       \tilde{Z}^2_{3 i} \Big{]}  \\
\hline (11)&\prodlim_{i=n+6}^{n+m+4} Z^2_{n+5,i} \cdot
       \prodlim_{n+6 \leq i < j \leq n+m+4}^{} {(Z^2_{i j})}^{Z^{-2}_{5 i} Z^2_{5,n+5}}  \\
\hline (12)&\prodlim_{i=n+5}^{n+m+4} {({Z}^2_{1 i})}^{Z^{2}_{1,n+5} Z^2_{1 3}}
       \cdot \prodlim_{i=6}^{n+4} \Big{[}\prodlim_{j=n+5}^{n+m+4} {({Z}^2_{i j})}^{(\prodlim_{k=n+5}^{j}
       Z^{-2}_{k j})} \Big{]}\\
\hline
\end{array}
$$
The skeletons of $\tilde{Z}_{4 5}$ and $\tilde{Z}^2_{3 i}$ appear
in Figure \ref{2figsnew}.

\end{proof}

\subsection{The corresponding fundamental groups}

In this section, we compute the fundamental group of the complement
of the general arrangement  $T_{n,m}$, $n, m \geq 1$.

Note that the arrangement $T_{1,1}$ appears already in
\cite{AmGaTe}, therefore we cite here the  presentation of its
fundamental group:

\begin{prop}
Let $T_{1,1}$ be the curve defined in Lemma \ref{t11}, see Figure \ref{thm3-fig}.

Then:
$$\pi_1(\C\P^2-T_{1,1}) \cong \langle x_1 \rangle \oplus
\langle x_2,x_3 \ |  \ (x_2 x_3)^2=(x_3 x_2)^2 \rangle.$$
\end{prop}

Now, we prove Theorem \ref{presentation_Tnm} by computing the
simplified presentation of the fundamental group of $T_{n,m}$:\\

\noindent
{\bf Proof of Theorem \ref{presentation_Tnm}:}\\
Using the van Kampen Theorem on $\Delta_{T_{n,m}}^2$ (computed in Proposition \ref{bmf-tnm}),
we get the following presentation for $\pi_1(\C\P^2-T_{n,m})$:

\noindent
Generators: $\{ x_1, x_2, x_3, x_4, x_5, x_6, x_7, \dots,  x_{n+4},  x_{n+5}, \dots, x_{n+m+4}\}$. \\
Relations:
\begin{tiny}
\begin{enumerate}
\item $x_{2} = x_{3}$

\item $x_3= x_{1}^{-1} x_3^{-1} x_4^{-1} x_5^{-1} \cdot x_{n+4} x_{n+3} \cdots x_3 x_2 x_3^{-1} \cdots
      x_{n+3}^{-1} x_{n+4}^{-1} \cdot x_{5} x_{4} x_{3} x_{1}$

\item $x_{4} x_{3} x_2 x_{3}^{-1} x_4 x_{3} x_{2}^{-1} x_{3}^{-1} x_{4}^{-1} = x_{n+5} x_{5} x_{n+5}^{-1}$

\item $x_3 x_{1} x_3^{-1} x_1^{-1} x_3^{-1} x_{4} x_{3} x_1 x_3 x_1^{-1} x_{3}^{-1} =
      x_{5}^{-1} \cdot x_{n+5}^{-1} x_{n+6}^{-1} \cdots  x_{n+m+4}^{-1} \cdot x_{n+5} x_{5} x_{n+5}^{-1}
      \cdot x_{n+m+4} \cdots  x_{n+6} x_{n+5} \cdot x_5$

\item $(x_{3} x_{2} x_3^{-1} x_4)^2 = (x_{4} x_{3} x_{2} x_3^{-1})^2$

\item $(x_{3} x_{1} x_3 x_1^{-1} x_3^{-1} x_4)^2 = (x_{4} x_{3} x_{1} x_3 x_1^{-1} x_3^{-1})^2$

\item $(x_{1} x_{3})^2 = (x_{3} x_{1})^2$

\item $(x_2 \cdot x_3^{-1} x_{4}^{-1} \cdots x_{i-1}^{-1} x_i x_{i-1} \cdots x_{4} x_3)^2 =
       (x_3^{-1} x_{4}^{-1} \cdots x_{i-1}^{-1} x_i x_{i-1} \cdots x_{4} x_3 \cdot x_{2})^2$,\\
       where $6 \leq i\leq n+4$

\item $(x_{5} x_{n+5})^2 = (x_{n+5} x_{5})^2$

\item $(x_{5} x_{n+5}^{-1} x_i x_{n+5})^2 = (x_{n+5}^{-1} x_{i} x_{n+5} x_5)^2$, where $n+6 \leq i \leq n+m+4$

\item $[x_{3} x_1 x_{3}^{-1}, x_{5}]=e$

\item $[x_{5} x_4 x_{5}^{-1}, x_i]=e$, where $6 \leq i \leq n+4$

\item $[x_{5} x_{3} x_{1} x_{3}^{-1} x_{5} x_{3} x_{1}^{-1} x_3^{-1} x_5^{-1}, x_i]=e$, where $6 \leq i \leq n+4$

\item $[x_5 x_{3} x_{1} x_{3}^{-1} x_{5}^{-1}, x_i]=e$, where $6 \leq i\leq n+4$.

\item $[x_{j-1} x_{j-2} \cdots x_4 x_3 x_2 x_{3}^{-1} x_{4}^{-1} \cdots  x_{i-2}^{-1} x_{i-1}^{-1} x_i
      x_{i-1} x_{i-2} \cdots  x_{4} x_{3} x_2^{-1}  x_3^{-1}  x_4^{-1}  \cdots x_{j-2}^{-1} x_{j-1}^{-1},
      x_j]=e$, where $6 \leq i < j \leq n+4$

\item $[x_4 x_{3} x_{2} x_{3}^{-1} x_{4}^{-1}, x_{n+5}]=e$

\item $[x_{3} x_{1} x_{3} x_{1}^{-1} x_{3}^{-1}, x_{n+5}]=e$

\item $[x_{5} x_{4} x_3 x_2 x_{3}^{-1} x_{4}^{-1} x_{5}^{-1}, x_{n+5}^{-1} x_{n+6}^{-1} \cdots  x_{i-2}^{-1}
      x_{i-1}^{-1} x_i x_{i-1} x_{i-2} \cdots x_{n+6} x_{n+5}]=e$, \\
      where $n+6 \leq i \leq n+m+4$

\item $[x_{5} x_{3} x_1 x_3 x_{1}^{-1} x_{3}^{-1} x_{5}^{-1}, x_{n+5}^{-1} x_{n+6}^{-1} \cdots  x_{i-2}^{-1}
      x_{i-1}^{-1} x_i x_{i-1} x_{i-2} \cdots x_{n+6} x_{n+5}]=e$, \\
      where $n+6 \leq i \leq n+m+4$

\item $[x_{n+5}, x_{i}]=e$, where $n+6 \leq i \leq n+m+4$

\item $[x_{n+5} x_{5}^{-1} x_{n+5}^{-1} x_{i} x_{n+5} x_{5} x_{n+5}^{-1}, x_j]=e$, where $n+6 \leq i < j \leq n+m+4$

\item $[x_{3} x_{1} x_{3}^{-1}, x_{n+5}^{-1} x_{i} x_{n+5}]=e$, where $n+5 \leq i \leq n+m+4$

\item $[x_{i}, x_{n+5}^{-1} x_{n+6}^{-1} \cdots x_{j-1}^{-1} x_{j} x_{j-1} \cdots x_{n+6} x_{n+5}]=e$,
      where $6 \leq i \leq n+4$ and \\
      $n+5 \leq j \leq n+m+4$

\item $x_{n+m+4} x_{n+m+3} \cdots x_4 x_3 x_2 x_1=e$. \mbox{(Projective relation)}.
\end{enumerate}
\end{tiny}

Since $x_2=x_3$ and $(x_{1} x_{3})^2 = (x_{3} x_{1})^2$ (by
Relations (1) and (7)),
we  derive a simpler presentation as follows (we combine Relations (5) and (7) into one set):

\noindent
Generators: $\{ x_1, x_2, x_4, x_5, x_6, x_7, \dots ,x_{n+4},  x_{n+5}, \dots, x_{n+m+4}\}$. \\
Relations:
\begin{tiny}
\begin{enumerate}
\item $x_1^{-1} x_2 x_1= x_4^{-1} x_5^{-1} \cdot x_{n+4} x_{n+3} \cdots x_4 x_2 x_4^{-1} \cdots
      x_{n+3}^{-1} x_{n+4}^{-1} \cdot x_{5} x_{4}$

\item $x_{4} x_2  x_4  x_{2}^{-1} x_{4}^{-1} = x_{n+5} x_{5} x_{n+5}^{-1}$

\item $x_1^{-1} x_2^{-1} x_1 x_{4} x_1^{-1} x_{2} x_{1} =
      x_{5}^{-1} \cdot x_{n+5}^{-1} x_{n+6}^{-1} \cdots  x_{n+m+4}^{-1} \cdot x_{n+5} x_{5} x_{n+5}^{-1}
      \cdot x_{n+m+4} \cdots  x_{n+6} x_{n+5} \cdot x_5$

\item $(x_{2} x_i)^2 = (x_{i} x_{2})^2$, where $i= 1, 4$

\item $(x_1^{-1} x_2 x_{1} x_4)^2 = (x_{4} x_1^{-1} x_2 x_{1})^2$

\item $(x_{4}^{-1} \cdots x_{i-1}^{-1} x_i x_{i-1} \cdots x_{4} x_2)^2 =
       (x_2^{-1} x_{4}^{-1} \cdots x_{i-1}^{-1} x_i x_{i-1} \cdots x_{4}
       x_{2}^2)^2$,\\
       where $6 \leq i\leq n+~4$

\item $(x_{5} x_{n+5})^2 = (x_{n+5} x_{5})^2$

\item $(x_{5} x_{n+5}^{-1} x_i x_{n+5})^2 = (x_{n+5}^{-1} x_{i} x_{n+5} x_5)^2$, where $n+6 \leq i \leq n+m+4$

\item $[x_{2} x_1 x_{2}^{-1}, x_{5}]=e$

\item $[x_{5} x_4 x_{5}^{-1}, x_i]=e$, where $6 \leq i \leq n+4$

\item $[x_{5} x_{2} x_{1} x_{2}^{-1} x_{5} x_{2} x_{1}^{-1} x_2^{-1} x_5^{-1}, x_i]=e$, where $6 \leq i \leq n+4$

\item $[x_5 x_{2} x_{1} x_{2}^{-1} x_{5}^{-1}, x_i]=e$, where $6 \leq i\leq n+4$.

\item $[x_{j-1} x_{j-2} \cdots x_4 x_2 x_{4}^{-1} \cdots  x_{i-2}^{-1} x_{i-1}^{-1} x_i
      x_{i-1} x_{i-2} \cdots  x_{4} x_2^{-1} x_4^{-1}  \cdots x_{j-2}^{-1} x_{j-1}^{-1},
      x_j]=e$, \\
      where $6 \leq i < j \leq n+4$

\item $[x_4 x_{2} x_{4}^{-1}, x_{n+5}]=e$

\item $[x_{1}^{-1} x_{2} x_{1}, x_{n+5}]=e$

\item $[x_{5} x_{4} x_2 x_{4}^{-1} x_{5}^{-1}, x_{n+5}^{-1}
x_{n+6}^{-1} \cdots  x_{i-2}^{-1} x_{i-1}^{-1} x_i x_{i-1} x_{i-2}
\cdots x_{n+6} x_{n+5}]=e$, \\
where $n+6 \leq i \leq n+m+4$

\item $[x_{5} x_{1}^{-1} x_{2} x_{1} x_{5}^{-1}, x_{n+5}^{-1} x_{n+6}^{-1} \cdots  x_{i-2}^{-1}
      x_{i-1}^{-1} x_i x_{i-1} x_{i-2} \cdots x_{n+6} x_{n+5}]=e$, \\
      where $n+6 \leq i \leq n+m+4$

\item $[x_{n+5}, x_{i}]=e$, where $n+6 \leq i \leq n+m+4$

\item $[x_{n+5} x_{5}^{-1} x_{n+5}^{-1} x_{i} x_{n+5} x_{5} x_{n+5}^{-1}, x_j]=e$, where $n+6 \leq i < j \leq n+m+4$

\item $[x_{2} x_{1} x_{2}^{-1}, x_{n+5}^{-1} x_{i} x_{n+5}]=e$, where $n+5 \leq i \leq n+m+4$

\item $[x_{i}, x_{n+5}^{-1} x_{n+6}^{-1} \cdots x_{j-1}^{-1} x_{j}
x_{j-1} \cdots x_{n+6} x_{n+5}]=e$,\\
      where $6 \leq i \leq n+4$ and $n+5 \leq j \leq n+m+4$

\item $x_{n+m+4} x_{n+m+3} \cdots x_{5} x_4 x_2^2 x_1=e$.
\end{enumerate}
\end{tiny}

Now we apply Relations (4) and (18) to get a much convenient
presentation (we also combine Relations (7) and (8) into one set):

\begin{tiny}
\begin{enumerate}
\item $x_1^{-1} x_2 x_1= x_4^{-1} x_5^{-1} \cdot x_{n+4} x_{n+3} \cdots x_4 x_2 x_4^{-1} \cdots
      x_{n+3}^{-1} x_{n+4}^{-1} \cdot x_{5} x_{4}$

\item $x_{2}^{-1} x_{4} x_{2} = x_{n+5} x_{5} x_{n+5}^{-1}$

\item $x_1^{-1} x_2^{-1} x_1 x_{4} x_1^{-1} x_{2} x_{1} =
      x_{5}^{-1} \cdot  x_{n+6}^{-1} \cdots  x_{n+m+4}^{-1} \cdot x_{5}
      \cdot x_{n+m+4} \cdots  x_{n+6} \cdot x_5$

\item $(x_{2} x_i)^2 = (x_{i} x_{2})^2$, where $i=1, 4$

\item $(x_1^{-1} x_2 x_{1} x_4)^2 = (x_{4} x_1^{-1} x_2 x_{1})^2$

\item $(x_{4}^{-1} \cdots x_{i-1}^{-1} x_i x_{i-1} \cdots x_{4} x_2)^2 =
       (x_2^{-1} x_{4}^{-1} \cdots x_{i-1}^{-1} x_i x_{i-1} \cdots x_{4}  x_{2}^2)^2$,\\
       where $6 \leq i\leq n+~4$

\item $(x_{5} x_i)^2 = (x_{i} x_5)^2$, where $n+5 \leq i \leq n+m+4$

\item $[x_{2} x_1 x_{2}^{-1}, x_{5}]=e$

\item $[x_{5} x_4 x_{5}^{-1}, x_i]=e$, where $6 \leq i \leq n+4$

\item $[x_{5} x_{2} x_{1} x_{2}^{-1} x_{5} x_{2} x_{1}^{-1} x_2^{-1} x_5^{-1}, x_i]=e$, where $6 \leq i \leq n+4$

\item $[x_5 x_{2} x_{1} x_{2}^{-1} x_{5}^{-1}, x_i]=e$, where $6 \leq i\leq n+4$.

\item $[x_{j-1} x_{j-2} \cdots x_4 x_2 x_{4}^{-1} \cdots  x_{i-2}^{-1} x_{i-1}^{-1} x_i
      x_{i-1} x_{i-2} \cdots  x_{4} x_2^{-1} x_4^{-1}  \cdots x_{j-2}^{-1} x_{j-1}^{-1},
      x_j]=e$,\\ where $6 \leq i < j \leq n+4$

\item $[x_4 x_{2} x_{4}^{-1}, x_{n+5}]=e$

\item $[x_{1}^{-1} x_{2} x_{1}, x_{n+5}]=e$

\item $[x_{5} x_{4} x_2 x_{4}^{-1} x_{5}^{-1}, x_{n+6}^{-1} \cdots  x_{i-2}^{-1}
x_{i-1}^{-1} x_i x_{i-1} x_{i-2} \cdots x_{n+6}]=e$, where $n+6 \leq i \leq n+m+4$

\item $[x_{5} x_{1}^{-1} x_{2} x_{1} x_{5}^{-1}, x_{n+6}^{-1} \cdots  x_{i-2}^{-1}
      x_{i-1}^{-1} x_i x_{i-1} x_{i-2} \cdots x_{n+6}]=e$, where $n+6 \leq i \leq n+m+4$

\item $[x_{n+5}, x_{i}]=e$, where $n+6 \leq i \leq n+m+4$

\item $[x_{5}^{-1} x_{i} x_{5}, x_j]=e$, where $n+6 \leq i < j \leq n+m+4$

\item $[x_{2} x_{1} x_{2}^{-1}, x_{i}]=e$, where $n+5 \leq i \leq n+m+4$

\item $[x_{i}, x_{n+6}^{-1} \cdots x_{j-1}^{-1} x_{j} x_{j-1} \cdots x_{n+6}]=e$,
      where $6 \leq i \leq n+4$ and $n+5 \leq j \leq n+m+4$

\item $x_{n+m+4} x_{n+m+3} \cdots x_{5} x_4 x_2^2 x_1=e$.
\end{enumerate}
\end{tiny}

By Relations (9) we get the following simplified relations (we
combine Relations (8), (11), (19) into one set):

\begin{tiny}
\begin{enumerate}
\item $x_1^{-1} x_2 x_1= x_4^{-1} x_5^{-1} \cdot x_{n+4} x_{n+3} \cdots x_4 x_2 x_4^{-1} \cdots
      x_{n+3}^{-1} x_{n+4}^{-1} \cdot x_{5} x_{4}$

\item $x_{2}^{-1} x_{4} x_{2} = x_{n+5} x_{5} x_{n+5}^{-1}$

\item $x_1^{-1} x_2^{-1} x_1 x_{4} x_1^{-1} x_{2} x_{1} =
      x_{5}^{-1} \cdot  x_{n+6}^{-1} \cdots  x_{n+m+4}^{-1} \cdot x_{5}
      \cdot x_{n+m+4} \cdots  x_{n+6} \cdot x_5$

\item $(x_{2} x_i)^2 = (x_{i} x_{2})^2$, where $i=1, 4$

\item $(x_1^{-1} x_2 x_{1} x_4)^2 = (x_{4} x_1^{-1} x_2 x_{1})^2$

\item $(x_{4}^{-1} \cdots x_{i-1}^{-1} x_i x_{i-1} \cdots x_{4} x_2)^2 =
       (x_2^{-1} x_{4}^{-1} \cdots x_{i-1}^{-1} x_i x_{i-1} \cdots x_{4}
       x_{2}^2)^2$,\\
       where $6 \leq i\leq n+~4$

\item $(x_{5} x_i)^2 = (x_{i} x_5)^2$, where $n+5 \leq i \leq n+m+4$

\item $[x_{2} x_1 x_{2}^{-1}, x_{i}]=e$, where $5 \leq i \leq n+m+4$

\item $[x_{5} x_4 x_{5}^{-1}, x_i]=e$, where $6 \leq i \leq n+4$

\item $[x_{5}, x_i]=e$, where $6 \leq i \leq n+4$

\item $[x_{j-1} x_{j-2} \cdots x_4 x_2 x_{4}^{-1} \cdots  x_{i-2}^{-1} x_{i-1}^{-1} x_i
      x_{i-1} x_{i-2} \cdots  x_{4} x_2^{-1} x_4^{-1}  \cdots x_{j-2}^{-1} x_{j-1}^{-1},
      x_j]=e$, \\
      where $6 \leq i < j \leq n+4$

\item $[x_4 x_{2} x_{4}^{-1}, x_{n+5}]=e$

\item $[x_{1}^{-1} x_{2} x_{1}, x_{n+5}]=e$

\item $[x_{5} x_{4} x_2 x_{4}^{-1} x_{5}^{-1}, x_{n+6}^{-1} \cdots  x_{i-2}^{-1}
x_{i-1}^{-1} x_i x_{i-1} x_{i-2} \cdots x_{n+6}]=e$, where $n+6 \leq i \leq n+m+4$

\item $[x_{5} x_{1}^{-1} x_{2} x_{1} x_{5}^{-1}, x_{n+6}^{-1} \cdots  x_{i-2}^{-1}
      x_{i-1}^{-1} x_i x_{i-1} x_{i-2} \cdots x_{n+6}]=e$, where $n+6 \leq i \leq n+m+4$

\item $[x_{n+5}, x_{i}]=e$, where $n+6 \leq i \leq n+m+4$

\item $[x_{5}^{-1} x_{i} x_{5}, x_j]=e$, where $n+6 \leq i < j \leq n+m+4$

\item $[x_{i}, x_{n+6}^{-1} \cdots x_{j-1}^{-1} x_{j} x_{j-1} \cdots x_{n+6}]=e$,
      where $6 \leq i \leq n+4$ and $n+5 \leq j \leq n+m+4$

\item $x_{n+m+4} x_{n+m+3} \cdots x_{5} x_4 x_2^2 x_1=e$.
\end{enumerate}
\end{tiny}

By Relations (10), Relations (9) become $[x_4, x_i]=e$, where $6 \leq i \leq n+4$.
These relations enable us to proceed in simplification:

\begin{tiny}
\begin{enumerate}
\item $x_1^{-1} x_2 x_1= x_{n+4} x_{n+3} \cdots x_6 x_2 x_6^{-1} \cdots
      x_{n+3}^{-1} x_{n+4}^{-1}$

\item $x_{2}^{-1} x_{4} x_{2} = x_{n+5} x_{5} x_{n+5}^{-1}$

\item $x_1^{-1} x_2^{-1} x_1 x_{4} x_1^{-1} x_{2} x_{1} =
      x_{5}^{-1} \cdot  x_{n+6}^{-1} \cdots  x_{n+m+4}^{-1} \cdot x_{5}
      \cdot x_{n+m+4} \cdots  x_{n+6} \cdot x_5$

\item $(x_{2} x_i)^2 = (x_{i} x_{2})^2$, where $i= 1,4$

\item $(x_1^{-1} x_2 x_{1} x_4)^2 = (x_{4} x_1^{-1} x_2 x_{1})^2$

\item $(x_{6}^{-1} \cdots x_{i-1}^{-1} x_i x_{i-1} \cdots x_{6} x_2)^2 =
       (x_2^{-1} x_{6}^{-1} \cdots x_{i-1}^{-1} x_i x_{i-1} \cdots x_{6}  x_{2}^2)^2$,\\
       where $6 \leq i\leq n+4$

\item $(x_{5} x_i)^2 = (x_{i} x_5)^2$, where $n+5 \leq i \leq n+m+4$

\item $[x_{2} x_1 x_{2}^{-1}, x_{i}]=e$, where $5 \leq i \leq n+m+4$

\item $[x_4, x_i]=e$, where $6 \leq i \leq n+4$

\item $[x_{5}, x_i]=e$, where $6 \leq i \leq n+4$

\item $[x_{j-1} x_{j-2} \cdots x_6 x_2 x_{6}^{-1} \cdots  x_{i-2}^{-1} x_{i-1}^{-1} x_i
      x_{i-1} x_{i-2} \cdots  x_{6} x_2^{-1} x_6^{-1}  \cdots x_{j-2}^{-1} x_{j-1}^{-1},
      x_j]=e$, \\
      where $6 \leq i < j \leq n+4$

\item $[x_4 x_{2} x_{4}^{-1}, x_{n+5}]=e$

\item $[x_{1}^{-1} x_{2} x_{1}, x_{n+5}]=e$

\item $[x_{5} x_{4} x_2 x_{4}^{-1} x_{5}^{-1}, x_{n+6}^{-1} \cdots  x_{i-2}^{-1}
x_{i-1}^{-1} x_i x_{i-1} x_{i-2} \cdots x_{n+6}]=e$, where $n+6 \leq i \leq n+m+4$

\item $[x_{5} x_{1}^{-1} x_{2} x_{1} x_{5}^{-1}, x_{n+6}^{-1} \cdots  x_{i-2}^{-1}
      x_{i-1}^{-1} x_i x_{i-1} x_{i-2} \cdots x_{n+6}]=e$, where $n+6 \leq i \leq n+m+4$

\item $[x_{n+5}, x_{i}]=e$, where $n+6 \leq i \leq n+m+4$

\item $[x_{5}^{-1} x_{i} x_{5}, x_j]=e$, where $n+6 \leq i < j \leq n+m+4$

\item $[x_{i}, x_{n+6}^{-1} \cdots x_{j-1}^{-1} x_{j} x_{j-1} \cdots x_{n+6}]=e$,
      where $6 \leq i \leq n+4$ and $n+5 \leq j \leq n+m+4$

\item $x_{n+m+4} x_{n+m+3} \cdots x_{5} x_4 x_2^2 x_1=e$.
\end{enumerate}
\end{tiny}

We treat Relations (18). Substituting $i=6$ and $j=n+5$, we
get $[x_{6}, x_{n+5}]=e$. Substituting $i=6$ and $j=n+6$, we get
$[x_{6}, x_{n+6}]=e$. If we take $i=6$ and $j=n+7$, we get
$[x_{6}, x_{n+6}^{-1} x_{n+7} x_{n+6}]=e$, which is $[x_{6},
x_{n+7}]=e$. Continuing this process by substituting $i=6$ and $
n+8\leq j \leq n+m+4$, we get $[x_{6}, x_{j}]=e$ for $ n+8 \leq j \leq
n+m+4$. Now we substitute $i=7$ and $ n+5 \leq j \leq n+m+4$. From
these substitutions we get:  $[x_{7}, x_{j}]=e$ for $ n+5 \leq j \leq
 n+m+4$. In the same manner, we get $[x_{i}, x_{j}]=e$ for
$8 \leq i \leq n+4$ and $n+5 \leq j \leq n+m+4$. This enables us
to simplify the above list of relations (we combine Relations (9),
(10) and (18) into one set):

\begin{tiny}
\begin{enumerate}
\item $x_1^{-1} x_2 x_1= x_{n+4} x_{n+3} \cdots x_6 x_2 x_6^{-1} \cdots
      x_{n+3}^{-1} x_{n+4}^{-1}$

\item $x_{2}^{-1} x_{4} x_{2} = x_{n+5} x_{5} x_{n+5}^{-1}$

\item $x_1^{-1} x_2^{-1} x_1 x_{4} x_1^{-1} x_{2} x_{1} =
      x_{5}^{-1} \cdot  x_{n+6}^{-1} \cdots  x_{n+m+4}^{-1} \cdot x_{5}
      \cdot x_{n+m+4} \cdots  x_{n+6} \cdot x_5$

\item $(x_{2} x_i)^2 = (x_{i} x_{2})^2$, where $i= 1,4$

\item $(x_1^{-1} x_2 x_{1} x_4)^2 = (x_{4} x_1^{-1} x_2 x_{1})^2$

\item $(x_{6}^{-1} \cdots x_{i-1}^{-1} x_i x_{i-1} \cdots x_{6} x_2)^2 =
       (x_2^{-1} x_{6}^{-1} \cdots x_{i-1}^{-1} x_i x_{i-1} \cdots x_{6}  x_{2}^2)^2$,\\
       where $6 \leq i\leq n+4$

\item $(x_{5} x_i)^2 = (x_{i} x_5)^2$, where $n+5 \leq i \leq n+m+4$

\item $[x_{2} x_1 x_{2}^{-1}, x_{i}]=e$, where $5 \leq i \leq n+m+4$

\item $[x_i, x_j]=e$, where $6 \leq i \leq n+4$ and $j=4, 5, n+5, \dots, n+m+4$

\item $[x_{j-1} x_{j-2} \cdots x_6 x_2 x_{6}^{-1} \cdots  x_{i-2}^{-1} x_{i-1}^{-1} x_i
      x_{i-1} x_{i-2} \cdots  x_{6} x_2^{-1} x_6^{-1}  \cdots x_{j-2}^{-1} x_{j-1}^{-1},
      x_j]=e$, \\ where $6 \leq i < j \leq n+4$

\item $[x_4 x_{2} x_{4}^{-1}, x_{n+5}]=e$

\item $[x_{1}^{-1} x_{2} x_{1}, x_{n+5}]=e$

\item $[x_{5} x_{4} x_2 x_{4}^{-1} x_{5}^{-1}, x_{n+6}^{-1} \cdots  x_{i-2}^{-1}
x_{i-1}^{-1} x_i x_{i-1} x_{i-2} \cdots x_{n+6}]=e$, where $n+6 \leq i \leq n+m+4$

\item $[x_{5} x_{1}^{-1} x_{2} x_{1} x_{5}^{-1}, x_{n+6}^{-1} \cdots  x_{i-2}^{-1}
      x_{i-1}^{-1} x_i x_{i-1} x_{i-2} \cdots x_{n+6}]=e$, where $n+6 \leq i \leq n+m+4$

\item $[x_{n+5}, x_{i}]=e$, where $n+6 \leq i \leq n+m+4$

\item $[x_{5}^{-1} x_{i} x_{5}, x_j]=e$, where $n+6 \leq i < j \leq n+m+4$

\item $x_{n+m+4} x_{n+m+3} \cdots x_{5} x_4 x_2^2 x_1=e$.
\end{enumerate}
\end{tiny}

Now we simplify Relations (6) and (10) (we have used this
trick of simplification in the computation of
$\pi_1(\C^2-C_n)$ in Section \ref{sec2}).

We start by substituting $i=6$ in Relations (6) to get $(x_2
x_6)^2=(x_6 x_2)^2$. We continue with $i=6, j=7$ in Relations (10)
to get $[x_2^{-1} x_6 x_2, x_7]=e$ (using $(x_2 x_6)^2=(x_6
x_2)^2$). Using these two relations, we simplify Relations (6) for
$i=7$: $(x_2 x_6^{-1} x_7 x_6)^2 = (x_6^{-1} x_7 x_6 x_2)^2$. We
add the expression $x_2^{-1} x_2=e$ in four locations as follows:
$$
x_6 x_2 x_6^{-1} (x_2^{-1} x_2) x_7 x_6 x_2 x_6^{-1} (x_2^{-1} x_2) x_7 = x_7 x_6 x_2 x_6^{-1}
(x_2^{-1} x_2) x_7 x_6 x_2 x_6^{-1} (x_2^{-1} x_2).
$$
Since $(x_2 x_6)^2=(x_6 x_2)^2$, we can rewrite the relation as:
$$
x_2^{-1} x_6^{-1} x_2 x_6 x_2 x_7 x_2^{-1} x_6^{-1} x_2 x_6  x_2 x_7 = x_7 x_2^{-1} x_6^{-1} x_2
x_6 x_2 x_7 x_2^{-1} x_6^{-1} x_2 x_6 x_2.
$$
We use the relation $[x_2^{-1} x_6 x_2, x_7]=e$ to get:
$$
x_2^{-1} x_6^{-1} x_2 x_6 x_2 x_2^{-1} x_6^{-1} x_2 x_7 x_6  x_2 x_7 = x_2^{-1} x_6^{-1} x_2
x_7 x_6 x_2 x_2^{-1} x_6^{-1} x_2 x_7 x_6 x_2,
$$
namely $x_2 x_7 x_6  x_2 x_7 = x_7 x_2 x_7 x_6 x_2$. We add the expression $x_2 x_2^{-1}=e$ to get
$x_2 x_7 (x_2 x_2^{-1}) x_6  x_2 x_7 = x_7 x_2 x_7 x_6 x_2$,
and using $[x_2^{-1} x_6 x_2, x_7]=e$, we get $(x_2 x_7)^2=(x_7 x_2)^2$.

Now we proceed with $i=6, j=8$ and then $i=7, j=8$ in Relations (10) to get
$[x_2^{-1} x_6 x_2, x_8]=e$ and $[x_2^{-1} x_7 x_2, x_8]=e$ respectively.
These relations enable us to simplify the relation $(x_2 x_6^{-1} x_7^{-1} x_8 x_7 x_6)^2=
(x_6^{-1} x_7^{-1} x_8 x_7 x_6 x_2)^2$ (which we get from Relations (6) for $i=8$) to
$(x_2 x_8)^2=(x_8 x_2)^2$.

In the same manner, we conclude that Relations (6) and (10) can be simplified to
$(x_2 x_i)^2=(x_i x_2)^2$ for $6 \leq i \leq n+4$ and $[x_2^{-1} x_i x_2, x_j]=e$ for
$6 \leq i < j \leq n+4$ respectively.

\smallskip
 Relation (12) gets the form $[x_2, x_{n+_5}]=e$ by
substituting Relation (1) in it and by using Relations (9).

Relations (14) can be rewritten by substituting Relation (1) in them (for $n+6 \leq i \leq n+m+4$):
$$
[x_{5} x_{n+4} x_{n+3} \cdots x_6 x_2 x_6^{-1} \cdots x_{n+3}^{-1} x_{n+4}^{-1} x_{5}^{-1},
x_{n+6}^{-1} \cdots  x_{i-2}^{-1} x_{i-1}^{-1} x_i x_{i-1} x_{i-2} \cdots x_{n+6}]=e.
$$
Now, by Relations (9), we get:
$$
[x_{5} x_2 x_{5}^{-1}, x_{n+6}^{-1} \cdots  x_{i-2}^{-1} x_{i-1}^{-1} x_i x_{i-1} x_{i-2} \cdots x_{n+6}]=e.
$$
For $i=n+6$, we get $[x_{5} x_2 x_{5}^{-1}, x_{n+6}]=e$. Now, for $i=n+7$, we get
$[x_{5} x_2 x_{5}^{-1}, x_{n+6}^{-1} x_{n+7} x_{n+6}]=e$. By  $[x_{5} x_2 x_{5}^{-1}, x_{n+6}]=e$,
it is simplified to $[x_{5} x_2 x_{5}^{-1}, x_{n+7}]=e$. In the same way, we get that Relations (14) are equivalent
to $[x_{5} x_2 x_{5}^{-1}, x_{i}]=e$, for  $n+6 \leq i \leq n+m+4$.\\

Therefore we have the following list of relations:
\begin{tiny}
\begin{enumerate}
\item $x_1^{-1} x_2 x_1= x_{n+4} x_{n+3} \cdots x_6 x_2 x_6^{-1} \cdots
      x_{n+3}^{-1} x_{n+4}^{-1}$

\item $x_{4}  = x_2 x_{n+5} x_{5} x_{n+5}^{-1} x_2^{-1}$

\item $x_1^{-1} x_2^{-1} x_1 x_{4} x_1^{-1} x_{2} x_{1} =
      x_{5}^{-1} \cdot  x_{n+6}^{-1} \cdots  x_{n+m+4}^{-1} \cdot x_{5}
      \cdot x_{n+m+4} \cdots  x_{n+6} \cdot x_5$

\item $(x_{2} x_i)^2 = (x_{i} x_{2})^2$, where $i= 1,4$

\item $(x_1^{-1} x_2 x_{1} x_4)^2 = (x_{4} x_1^{-1} x_2 x_{1})^2$

\item $(x_{2} x_i)^2 = (x_i x_2)^2$, where $6 \leq i\leq n+4$

\item $(x_{5} x_i)^2 = (x_{i} x_5)^2$, where $n+5 \leq i \leq n+m+4$

\item $[x_{2} x_1 x_{2}^{-1}, x_{i}]=e$, where $5 \leq i \leq n+m+4$

\item $[x_i, x_j]=e$, where $6 \leq i \leq n+4$ and $j=4, 5, n+5, \dots, n+m+4$

\item $[x_{2}^{-1} x_i x_{2}, x_{j}]=e$, where $6 \leq i < j \leq n+4$

\item $[x_4 x_{2} x_{4}^{-1}, x_{n+5}]=e$

\item $[x_{2}, x_{n+5}]=e$

\item $[x_{5} x_{4} x_2 x_{4}^{-1} x_{5}^{-1}, x_{n+6}^{-1} \cdots  x_{i-2}^{-1}
x_{i-1}^{-1} x_i x_{i-1} x_{i-2} \cdots x_{n+6}]=e$, where $n+6 \leq i \leq n+m+4$

\item $[x_{5} x_{2} x_{5}^{-1}, x_i]=e$, where $n+6 \leq i \leq n+m+4$

\item $[x_{n+5}, x_{i}]=e$, where $n+6 \leq i \leq n+m+4$

\item $[x_{5}^{-1} x_{i} x_{5}, x_j]=e$, where $n+6 \leq i < j \leq n+m+4$

\item $x_{n+m+4} x_{n+m+3} \cdots x_{5} x_4 x_2^2 x_1=e$.
\end{enumerate}
\end{tiny}

By Relation (2), we can omit $x_4$ and replace it everywhere by 
$$x_2 x_{n+5} x_5 x_{n+5}^{-1} x_2^{-1}.$$

We start with  Relations (4). Take $(x_{2} x_4)^2 = (x_{4}
x_{2})^2$. This relation is rewritten as:
\begin{eqnarray*}
x_2 \cdot x_2 x_{n+5} x_5 x_{n+5}^{-1} x_2^{-1} \cdot x_2 \cdot x_2 x_{n+5} x_5 x_{n+5}^{-1} x_2^{-1}=\\
=x_2 x_{n+5} x_5 x_{n+5}^{-1} x_2^{-1} \cdot x_2 \cdot x_2 x_{n+5} x_5 x_{n+5}^{-1} x_2^{-1} \cdot x_2.
\end{eqnarray*}
By Relation (12) we easily get:$(x_{2} x_5)^2 = (x_{5} x_{2})^2$.

Now we consider Relation (11):
\begin{eqnarray*}
e&=&[x_4 x_{2} x_{4}^{-1}, x_{n+5}]=\\&=&[x_2 x_{n+5} x_5 x_{n+5}^{-1}
x_2^{-1} \cdot x_2 \cdot x_2 x_{n+5} x_5^{-1} x_{n+5} x_2,
x_{n+5}]=\\ &=&[x_5 x_2 x_5^{-1},  x_{n+5}].
\end{eqnarray*}

By the above resulting relations, we can simplify also Relations (13):

\begin{tiny}
\begin{eqnarray*}
e&=&[x_{5} \cdot  x_2 x_{n+5} x_5 x_{n+5}^{-1} x_2^{-1} \cdot  x_2 \cdot x_2 x_{n+5} x_5^{-1}
x_{n+5}^{-1} x_2^{-1} \cdot x_{5}^{-1}, x_{n+6}^{-1} \cdots x_{i-1}^{-1} x_i x_{i-1} \cdots x_{n+6}]=\\
&=&[x_{2}, x_{n+6}^{-1} \cdots x_{i-1}^{-1} x_i x_{i-1} \cdots x_{n+6}],
\end{eqnarray*}
\end{tiny}

\noindent
for $n+6 \leq i \leq n+m+4$.
For $i=n+6$, we get $[x_2, x_{n+6}]=e$. For $i=n+7$, we get $[x_2, x_{n+6}^{-1} x_{n+7} x_{n+6}]=e$,
and by the previous relation, we get $[x_2, x_{n+7}]=e$. In the same manner, we get
 $[x_2, x_{i}]=e$, for  $n+6 \leq i \leq n+m+4$.

From Relations (9), we have $[x_i, x_4]=e$, where $6 \leq i \leq
n+4$. By Relation (2),  $e=[x_i, x_4]=[x_i, x_2 x_{n+5} x_{5}
x_{n+5}^{-1} x_2^{-1}]$. By Relations (9) and (12),  we get for $6
\leq i \leq n+4$:  $[x_i, x_2 x_{5} x_2^{-1}]=e$.

We rewrite now Relation (5), by substituting Relations (1) and (2) in it:

\begin{tiny}
\begin{eqnarray*}
x_{n+4}  \cdots x_6 x_2 x_6^{-1} \cdots  x_{n+4}^{-1} \cdot x_2 x_{n+5} x_{5} x_{n+5}^{-1} x_2^{-1} \cdot
x_{n+4}  \cdots x_6 x_2 x_6^{-1} \cdots  x_{n+4}^{-1} \cdot x_2 x_{n+5} x_{5} x_{n+5}^{-1} x_2^{-1}= \\
=x_2 x_{n+5} x_{5} x_{n+5}^{-1} x_2^{-1} \cdot x_{n+4}  \cdots x_6 x_2 x_6^{-1} \cdots  x_{n+4}^{-1} \cdot
x_2 x_{n+5} x_{5} x_{n+5}^{-1} x_2^{-1} \cdot x_{n+4}  \cdots x_6 x_2 x_6^{-1} \cdots  x_{n+4}^{-1}.
\end{eqnarray*}
\end{tiny}

Since $[x_{n+5}, x_{i}]=e$, for  $i= 2, 6, \dots, n+4$ (by Relations (9) and (12)),
\begin{eqnarray*}
x_{n+4}  \cdots x_6 x_2 x_6^{-1} \cdots  x_{n+4}^{-1}  x_2 x_{5}  x_2^{-1} \cdot
x_{n+4}  \cdots x_6 x_2 x_6^{-1} \cdots  x_{n+4}^{-1}  x_2  x_{5} x_2^{-1} = \\
=x_2 x_{5} x_2^{-1} x_{n+4}  \cdots x_6 x_2 x_6^{-1} \cdots  x_{n+4}^{-1} \cdot
x_2 x_{5} x_2^{-1} x_{n+4}  \cdots x_6 x_2 x_6^{-1} \cdots  x_{n+4}^{-1}.
\end{eqnarray*}
Now, since we proved above that $[x_i, x_2 x_{5} x_2^{-1}]$, for $6 \leq i \leq n+4$, we get
$x_2 x_{5} x_2  x_{5} =  x_{5}  x_2 x_{5} x_2$, which is a consequence from Relations (4) (see above).

\medskip
Therefore, we have:\\
Generators: $\{ x_1, x_2, x_5, x_6, x_7, \dots,  x_{n+4},  x_{n+5}, \dots, x_{n+m+4}\}$. \\
Relations:

\begin{tiny}
\begin{enumerate}
\item $x_1^{-1} x_2 x_1= x_{n+4} x_{n+3} \cdots x_6 x_2 x_6^{-1} \cdots x_{n+3}^{-1} x_{n+4}^{-1}$

\item $x_1^{-1} x_2^{-1} x_1 \cdot x_{2} x_{n+5} x_5 x_{n+5}^{-1} x_2^{-1} \cdot x_1^{-1} x_{2} x_{1} =
x_{5}^{-1} \cdot  x_{n+6}^{-1} \cdots  x_{n+m+4}^{-1} \cdot x_{5} \cdot x_{n+m+4} \cdots  x_{n+6} \cdot x_5$

\item $(x_{2} x_i)^2 = (x_{i} x_{2})^2$, where $i= 1,5, 6, \dots, n+4$

\item $(x_{5} x_i)^2 = (x_{i} x_5)^2$, where $n+5 \leq i \leq n+m+4$

\item $[x_{2} x_1 x_{2}^{-1}, x_{i}]=e$, where $5 \leq i \leq n+m+4$

\item $[x_i, x_j]=e$, where $6 \leq i \leq n+4$ and $j=5, n+5, \dots, n+m+4$

\item $[x_i, x_2 x_{5} x_2^{-1}]$, where $6 \leq i \leq n+4$

\item $[x_{2}^{-1} x_i x_{2}, x_{j}]=e$, where $6 \leq i < j \leq n+4$

\item $[x_2, x_{i}]=e$, where  $n+6 \leq i \leq n+m+4$

\item $[x_{5} x_{2} x_{5}^{-1}, x_i]=e$, where $n+5 \leq i \leq n+m+4$

\item $[x_{n+5}, x_{i}]=e$, where $i= 2, n+6, \dots, n+m+4$

\item $[x_{5}^{-1} x_{i} x_{5}, x_j]=e$, where $n+6 \leq i < j \leq n+m+4$

\item $x_{n+m+4} x_{n+m+3} \cdots x_{5} x_2 x_{n+5} x_{5} x_{n+5}^{-1} x_2 x_1=e$.
\end{enumerate}
\end{tiny}

Now we use the projective relation (Relation (13)) in order to
omit the generator $x_1$.

We rewrite  Relation (1), by substituting $x_1$:
\begin{eqnarray*}
& x_{n+m+4} \cdots x_{5} x_2 x_{n+5} x_{5} x_{n+5}^{-1} x_2 x_2
x_2^{-1}  x_{n+5} x_5^{-1}  x_{n+5}^{-1} x_2^{-1} x_5^{-1} \cdots  x_{n+m+4}^{-1}  = \\
& = x_{n+4}  \cdots x_6 x_2 x_6^{-1} \cdots  x_{n+4}^{-1}.
\end{eqnarray*}
By Relations (6) and (11), we get:
\begin{eqnarray*}
& x_{n+m+4} \cdots x_{n+5} x_{5} x_2 x_{n+5} x_{5}  x_2  x_5^{-1} x_{n+5}^{-1} x_2^{-1} x_5^{-1} x_{n+5}^{-1}
\cdots  x_{n+m+4}^{-1} = x_2 .
\end{eqnarray*}
By Relations (9) and (10), we have:
\begin{eqnarray*}
x_{n+5} x_{5} x_2 x_{5}  x_2  x_5^{-1} x_2^{-1} x_5^{-1} x_{n+5}^{-1} = x_2.
\end{eqnarray*}
Now, by Relations (3),
\begin{eqnarray*}
x_{n+5} x_2  x_{n+5}^{-1} =  x_2,
\end{eqnarray*}
which is known already in (11).\\

Now we simplify $(x_1 x_2)^2=(x_2 x_1)^2$, which appears in
Relations (3). By substituting $x_1$ and some immediate
cancellations,
\begin{eqnarray*}
x_{n+4} \cdots x_{6} x_5 x_{2} x_{n+5} x_5  x_{n+m+4} \cdots x_{n+6} x_{n+4} \cdots x_{6} = \\
=x_2^{-1} x_{n+4} \cdots x_{6} x_5 x_2 x_{n+5} x_{5}  x_{n+m+4} \cdots x_{n+6} x_{n+4} \cdots x_{6} x_2.
\end{eqnarray*}
By Relations (6) and (9),
$$x_{n+4} \cdots x_{6} x_5 x_{2} x_{n+5} x_5  x_{n+4} \cdots x_{6} =
x_2^{-1} x_{n+4} \cdots x_{6} x_5 x_2 x_{n+5} x_{5} x_{n+4} \cdots x_{6} x_2.$$
Since $[x_i, x_5]=e$ for  $6 \leq i \leq n+4$ (by Relations (6)), we get
$$x_{n+4} \cdots x_{6} x_5 x_{2} x_{n+5} x_{n+4} \cdots x_{6} x_5 =
x_2^{-1} x_{n+4} \cdots x_{6} x_5 x_2 x_{n+5} x_{n+4} \cdots x_{6} x_5 x_2,$$
which is:
$$x_{n+4} \cdots x_{6} x_5 x_{2} x_{n+5} x_{n+4} \cdots x_{6} x_5 x_2^{-1} x_5^{-1}=
x_2^{-1} x_{n+4} \cdots x_{6} x_5 x_2 x_{n+5} x_{n+4} \cdots x_{6}.$$
By Relations (6) and (10),
$$x_{n+4} \cdots x_{6} x_5 x_{2} x_{n+4} \cdots x_{6} x_5 x_2^{-1} x_5^{-1}=
x_2^{-1} x_{n+4} \cdots x_{6} x_5 x_2 x_{n+4} \cdots x_{6}.$$
Now we add $x_2 x_2^{-1} =e$ in two locations:
$$x_2 x_{n+4} \cdots x_{6} (x_2 x_2^{-1}) x_5 x_{2} x_{n+4} \cdots x_{6} x_5 x_2^{-1} x_5^{-1}=
x_{n+4} \cdots x_{6} (x_2 x_2^{-1}) x_5 x_2 x_{n+4} \cdots x_{6}.$$
By Relations (7),
$$x_2 x_{n+4} \cdots x_{6} x_2 x_{n+4} \cdots x_{6} x_2^{-1} x_5 x_2 x_5 x_2^{-1} x_5^{-1}=
x_{n+4} \cdots x_{6} x_2 x_{n+4} \cdots x_{6} x_2^{-1} x_5 x_2.$$
Since $(x_2 x_5)^2 = (x_5 x_2)^2$ (see Relations (3)), we get:
$$x _2 x_{n+4} \cdots x_{6} x_2 x_{n+4} \cdots x_{6}=
x_{n+4} \cdots x_{6} x_2 x_{n+4} \cdots x_{6} x_2.$$
Now we add again $x_2 x_2^{-1}=e$ as follows:
\begin{eqnarray*}
x_2 (x_2 x_2^{-1}) x_{n+4} (x_2 x_2^{-1}) \cdots (x_2 x_2^{-1}) x_7 (x_2 x_2^{-1}) x_{6} x_2 x_{n+4} \cdots x_{6}= \\
=x_{n+4} (x_2 x_2^{-1}) \cdots (x_2 x_2^{-1}) x_{6} x_2 x_{n+4} \cdots x_{6} x_2.
\end{eqnarray*}
We can rewrite the relation as follows:
\begin{eqnarray*}
x_2 (x_2^{-1} x_{n+4} x_2) x_2^{-1} \cdots x_2 (x_2^{-1} x_7 x_2) (x_2^{-1} x_{6} x_2) x_{n+4} \cdots x_{6}= \\
=(x_2^{-1} x_{n+4} x_2) x_2^{-1} \cdots x_2 (x_2^{-1} x_7 x_2) (x_2^{-1} x_{6} x_2) x_{n+4} \cdots x_{6} x_2,
\end{eqnarray*}
and using Relations (8), we get:
\begin{eqnarray*}
x_2 (x_2^{-1} x_{n+4} x_2) x_{n+4} \cdots (x_2^{-1} x_7 x_2)  x_7 (x_2^{-1} x_{6} x_2)  x_{6}=\\
=(x_2^{-1} x_{n+4} x_2) x_{n+4} \cdots (x_2^{-1} x_7 x_2) x_7 (x_2^{-1} x_{6} x_2) x_{6} x_2.
\end{eqnarray*}
Since by Relations (3) we have $(x_2 x_i)^2=(x_i x_2)^2$ for 
$6 \leq i \leq n+4$, this relation is redundant.

Relation (2) can be simplified by substituting Relation (1) and using Relations (6) and (11):

\begin{tiny}
\begin{eqnarray*}
x_{n+5} \cdot x_{n+4} \cdots x_6 \cdot x_{2}^{-1} \cdot x_6^{-1}  \cdots x_{n+4}^{-1} \cdot x_2  x_5 x_2^{-1}
\cdot  x_{n+4} \cdots x_6 \cdot x_{2} \cdot x_6^{-1}  \cdots x_{n+4}^{-1} \cdot x_{n+5}^{-1} =\\
=x_{5}^{-1} \cdot  x_{n+6}^{-1} \cdots  x_{n+m+4}^{-1} \cdot x_{5} \cdot x_{n+m+4} \cdots  x_{n+6} \cdot x_5.
\end{eqnarray*}
\end{tiny}
By Relations (6) and (7),
\begin{eqnarray*}
x_5 x_{n+5} x_5 x_{n+5}^{-1} x_5^{-1} = x_{n+6}^{-1} \cdots  x_{n+m+4}^{-1} \cdot x_{5} \cdot x_{n+m+4} \cdots  x_{n+6}.
\end{eqnarray*}
Using again  $(x_5 x_{n+5})^2 = (x_{n+5} x_5)^2$, we get:
\begin{equation*}
x_{n+5}^{-1} x_5 x_{n+5} = x_{n+6}^{-1} \cdots  x_{n+m+4}^{-1} \cdot x_{5} \cdot x_{n+m+4} \cdots  x_{n+6},
\end{equation*}
which can be rewritten as:
\begin{equation*}\label{3new}
{(*)} \ \ \ \ \ \ \ \ \ x_{n+m+4} \cdots  x_{n+6} \cdot x_5 \cdot x_{n+5} = x_{n+5} \cdot x_{5} \cdot x_{n+m+4} \cdots  x_{n+6}.
\end{equation*}

\medskip

We will see now that Relations (5) are all redundant. We split it 
into four subcases, according to the value of $i$:
\begin{itemize}
\item
If we substitute $i=5$ in Relations (5), we have $[x_1, x_2^{-1} x_5 x_2]=e$, which can be simplified to:
$$x_{n+m+4} \cdots x_{5} x_2 x_{n+5} x_{5} x_{n+5}^{-1} x_5 = x_2^{-1} x_5 x_2
x_{n+m+4} \cdots x_{5} x_2 x_{n+5} x_{5} x_{n+5}^{-1}.$$
By adding $x_5^{-1} x_5=e$ and using $(x_5 x_{n+5})^2 = (x_{n+5} x_5)^2$, we get:
\begin{eqnarray*}
x_{n+m+4} \cdots x_{n+6} x_{n+5} x_{n+4} \cdots x_6 x_{5} x_2 x_5^{-1} x_{n+5}^{-1} x_5 x_{n+5} x_5 = \\
=x_2^{-1} x_5 x_2 x_{n+m+4} \cdots x_{n+6} x_{n+5} x_{n+4} \cdots x_6 x_{5} x_2 x_5^{-1} x_{n+5}^{-1} x_{5} x_{n+5}.
\end{eqnarray*}
By Relations (6), (7) and (10),
\begin{eqnarray*}
x_{n+m+4} \cdots x_{n+6} x_{5} x_2 x_{n+5} x_5 x_{n+5}^{-1} x_5^{-1}=
x_2^{-1} x_5 x_2 x_{n+m+4} \cdots x_{n+6}  x_{5} x_2 x_5^{-1}.
\end{eqnarray*}
By Relations (10), (11) and $(x_2 x_5)^2=(x_5 x_2)^2$:
\begin{eqnarray*}
x_{n+m+4} \cdots x_{n+6} x_{5}  x_{n+5} x_5 x_{n+5}^{-1} x_5^{-1} x_2 = x_5 x_{2} x_{n+m+4} \cdots x_{n+6}.
\end{eqnarray*}
By Relations (9) and $(x_5 x_{n+5})^2 = (x_{n+5} x_5)^2$,
\begin{eqnarray*}
x_{n+5}^{-1} x_5 x_{n+5} =  x_{n+6}^{-1} \cdots x_{n+m+4}^{-1} x_5 x_{n+m+4} \cdots x_{n+6},
\end{eqnarray*}
which is known already (by $(*)$).

\item
For $6 \leq i \leq n+4$, we have $[x_1, x_2^{-1} x_i x_2]=e$.
By some cancellations, we have:
\begin{eqnarray*}
x_{n+m+4}\cdots x_{5} x_2 x_{n+5} x_{5} x_{n+5}^{-1} x_i =
x_2^{-1} x_i x_2 x_{n+m+4} \cdots x_{5} x_2 x_{n+5} x_{5} x_{n+5}^{-1}.
\end{eqnarray*}
By Relations (6) and (9), $x_{n+4}\cdots x_6 x_{5} x_2 x_i = x_2^{-1} x_i x_2 x_{n+4} \cdots x_6 x_{5} x_2$.
In this relation we substitute first $i=6$ to get (by Relations (8)):
$$x_6 x_{5} x_2 x_6 = x_2^{-1} x_6 x_2  x_6 x_{5} x_2,$$
which is $[x_6, x_2^{-1} x_5 x_2]=e$. This relation appears already in Relations (7).
Now we substitute $i=7$ to get (again by Relations (8)):
$$x_7 x_6 x_{5} x_2 x_7 = x_2^{-1} x_7 x_2 x_7 x_6 x_{5} x_2,$$
which is $x_6 x_{5} x_2 x_7 = x_2 x_7 x_2^{-1} x_6 x_{5} x_2$.
The addition of $x_2 x_2^{-1}=e$ in two locations, enables us to translate
$$x_2^{-1} x_6 (x_2 x_2^{-1}) x_{5} x_2 x_7 = x_7 x_2^{-1} x_6 (x_2 x_2^{-1}) x_{5} x_2$$
to $(x_2^{-1} x_6 x_2) (x_2^{-1} x_{5} x_2) x_7 = x_7 (x_2^{-1} x_6 x_2) (x_2^{-1} x_{5} x_2)$,
which is of course redundant by Relations (8). As for the other indices, all relations are redundant in a similar way.

\item
If $i=n+5$, we have by Relations (11): $[x_1 , x_{n+5}]=e$.
By some immediate cancellations, we get:
$$x_{n+m+4} x_{n+m+3} \cdots x_{5} x_2  =
x_{n+5} x_{n+m+4} x_{n+m+3} \cdots x_{5} x_2 x_{n+5} x_{5} x_{n+5}^{-1} x_5^{-1} x_{n+5}^{-1}.$$
By Relations (6) and (11), we have $x_{5} x_2  = x_{n+5} x_{5} x_2 x_{n+5} x_{5} x_{n+5}^{-1} x_5^{-1} x_{n+5}^{-1}$.
By Relations (4) and (10), the relation is redundant.

\item
For $n+6 \leq i \leq n+m+4$, we have by Relations (9): $[x_1, x_i]=e$. Substituting Relation (13), we have:
\begin{eqnarray*}
x_{n+m+4} \cdots x_{5} x_2 x_{n+5} x_{5} x_{n+5}^{-1} x_2 x_i =
x_i  x_{n+m+4} \cdots x_{5} x_2 x_{n+5} x_{5} x_{n+5}^{-1} x_2,
\end{eqnarray*}\\
which is (by Relations (6), (9) and (11)):
\begin{eqnarray*}
x_{n+m+4} \cdots x_{n+6} x_{5} x_2 x_{n+5} x_{5}  x_i =x_i  x_{n+m+4} \cdots x_{n+6} x_{5} x_2 x_{n+5} x_{5}.
\end{eqnarray*}
We add $x_5^{-1} x_5=e$ in two locations as follows:
\begin{eqnarray*}
x_{n+m+4} \cdots x_{n+6} x_{5} x_2 (x_5^{-1} x_5) x_{n+5} x_{5}  x_i =
x_i  x_{n+m+4} \cdots x_{n+6} x_{5} x_2 (x_5^{-1} x_5) x_{n+5} x_{5}.
\end{eqnarray*}
This enables us to use Relations (10) and to cancel $x_{5} x_2 x_5^{-1}$:
\begin{eqnarray*}
x_{n+m+4} \cdots x_{n+6} x_5 x_{n+5} x_{5}  x_i =x_i  x_{n+m+4} \cdots x_{n+6} x_5 x_{n+5} x_{5}.
\end{eqnarray*}
By Relation ($*$), we can rewrite it as:
\begin{eqnarray*}
x_{n+5} x_{5} x_{n+m+4} \cdots x_{n+6} x_5 x_i =x_i  x_{n+5} x_{5} x_{n+m+4} \cdots x_{n+6} x_5.
\end{eqnarray*}
By relations (11), we get  $[x_5^{-1} x_{i} x_{5}, x_{n+m+4} \cdots x_{n+6} x_5^2]=e$.

Now we substitute $i=n+6$:
\begin{eqnarray*}
x_5^{-1} x_{n+6} x_{5} x_{n+m+4} \cdots x_{n+6} x_5^2=
x_{n+m+4} \cdots x_{n+6} x_5^2 x_5^{-1} x_{n+6} x_{5}.
\end{eqnarray*}
By Relations (12), we get $x_{n+6} x_{5} x_{n+6} x_5= x_5 x_{n+6} x_5 x_{n+6}$,
which is known already (by Relations (4)).\\
Now we substitute $i=n+7$:
\begin{eqnarray*}
x_5^{-1} x_{n+7} x_{5} x_{n+m+4} \cdots x_{n+6} x_5=x_{n+m+4} \cdots x_{n+6} x_5 x_{n+7},
\end{eqnarray*}
Again, by Relations (12), we get:
\begin{eqnarray*}
x_5^{-1} x_{n+7} x_{5} x_{n+7} x_{n+6} x_5=x_{n+7} x_{n+6} x_5 x_{n+7},
\end{eqnarray*}
and by $(x_5 x_{n+7})^2 =(x_{n+7} x_5)^2$ (Relations (4)), we get:
\begin{eqnarray*}
x_5 x_{n+7} x_{5}^{-1} x_{n+6} x_5=x_{n+6} x_5 x_{n+7}.
\end{eqnarray*}
Now, since $[x_{n+7}, x_{5}^{-1} x_{n+6} x_5]=e$ (Relations (12)), this relation is redundant.

In a similar way, the relations $[x_1, x_i]=e$ for $n+6 \leq i \leq n+m+4$, are redundant.
\end{itemize}

\medskip

Therefore the resulting presentation is as follows:\\
Generators: $\{x_2, x_5, x_6, x_7, \dots ,x_{n+4},  x_{n+4},  x_{n+5}, \dots, x_{n+m+4}\}$. \\
Relations:

\begin{enumerate}

\item $x_{n+5}^{-1} x_5 x_{n+5} = x_{n+6}^{-1} \cdots  x_{n+m+4}^{-1} \cdot x_{5} \cdot x_{n+m+4} \cdots  x_{n+6}$

\item $(x_{2} x_i)^2 = (x_{i} x_{2})^2$, where $i=5, 6, \dots, n+4$

\item $(x_{5} x_i)^2 = (x_{i} x_5)^2$, where $n+5 \leq i \leq n+m+4$

\item $[x_i, x_j]=e$, where $6 \leq i \leq n+4$ and $j=5, n+5, \dots, n+m+4$

\item $[x_i, x_2 x_{5} x_2^{-1}]$, where $6 \leq i \leq n+4$

\item $[x_{2}^{-1} x_i x_{2}, x_{j}]=e$, where $6 \leq i < j \leq n+4$

\item $[x_2, x_{i}]=e$, where  $n+6 \leq i \leq n+m+4$

\item $[x_{5} x_{2} x_{5}^{-1}, x_i]=e$, where $n+5 \leq i \leq n+m+4$

\item $[x_{n+5}, x_{i}]=e$, where $i= 2, n+6, \dots, n+m+4$

\item $[x_{5}^{-1} x_{i} x_{5}, x_j]=e$, where $n+6 \leq i < j \leq n+m+4$,
\end{enumerate}

as stated in Theorem \ref{presentation_Tnm}. \hfill $\qed$

\medskip

As we have seen in Corollary \ref{big2}, the projective
fundamental groups of this section are big.

\begin{rem}
Note that if we substitute  $m=0$ in the presentation of Theorem \ref{presentation_Tnm},
we get exactly the presentation of $T_{n,0}$ given in Theorem \ref{presentation_Tn0}.
\end{rem}

\section{Acknowledgements}
We would like to thank Mutsuo Oka for correcting a wrong
assumption in the first version of this paper.

\end{document}